\documentclass[11pt, reqno]{amsart}
 \usepackage{amsmath}
 \usepackage{amssymb}
\usepackage{bm}

\DeclareMathAccent{\mathring}{\mathalpha}{operators}{"17}

 \usepackage{color}

\newcommand{\mysection}[1]{\section{#1}
      \setcounter{equation}{0}}

\newtheorem{theorem}{Theorem}[section]
\newtheorem{lemma}[theorem]{Lemma}

\newtheorem{corollary}[theorem]{Corollary} 

\theoremstyle{definition}
\newtheorem{assumption}{Assumption}[section]

\theoremstyle{remark}
\newtheorem{remark}{Remark}[section]

\newtheorem{example}{Example}[section]

\newcommand{\tr}{\text{\rm tr}\,}
\newcommand{\loc}{\text{\rm loc}}

 \makeatletter 
 \def\dashint{%  
 \operatorname%
 {\,\,\text{\bf--}\kern-.98em\DOTSI\intop\ilimits@\!\!}}
\def\ninf{\qopname\relax\@empty{inf\phantom{p}\!\!\!}}
\def\nelambda{\dot\lambda}
\def\melambda{\ddot\lambda}

 \makeatother

\newcommand\bbeta{\text{\raise-.2ex\hbox{$\bm{\beta}$}}}

\newcommand\bR{\mathbb{R}}

\newcommand\bQ{\mathbb{Q}}
\newcommand\bS{\mathbb{S}}

\newcommand\cF{\mathcal{F}}

\def\sft{{\sf t}}

\newcommand\uR{\underline{R}}

\newcommand\ulambda{\underline\lambda}

\newcommand\dist{{\rm dist}\,}

\begin{document}

\title[Potentials of It\^o's processes with drift in $L_{d+1}$]
{On potentials of It\^o's processes with drift in $L_{d+1}$}

\author{N.V. Krylov}
 
\email{nkrylov@umn.edu}
\address{127 Vincent Hall, University of Minnesota,
 Minneapolis, MN, 55455}
 
\keywords{It\^o's equations with singular drift, Potentials
of diffusion
processes}

\subjclass[2010]{60H10, 60J60}

\begin{abstract}
This paper is a natural continuation of \cite{Kr_20_2},
where strong Markov processes are constructed
in time inhomogeneous setting with Borel measurable
uniformly bounded and uniformly nondegenerate diffusion
and drift in $L_{d+1}(\mathbb{R}^{d+1})$. Here we study
some properties of these processes such as
the probability to pass through narrow tubes,
 higher summability of 
Green's functions, and so on.
The results seem to be new even if the diffusion
is constant.
\end{abstract}

\maketitle

\mysection{Introduction}

Let $\bR^{d}$ be a  Euclidean space of points
$x=(x^{1},...,x^{d})$, $d\geq 2$. Fix some 
 $p_{0},q_{0} \in[1,\infty)$  
such that
\begin{equation}
                                                    \label{5.10.1}
  \frac{d}{p_{0}}+\frac{1}{q_{0}}= 1.
\end{equation}

It is proved in \cite{Kr_20_2} that
It\^o's stochastic equations of the form
\begin{equation}
                                                 \label{11.29.2}
x _{s}=x  +\int_{0}^{s}\sigma (t +r,x_{r})\,dw_{r}
+\int_{0}^{s}b (t +r,x_{r}) \,dr
\end{equation}
admit  weak solutions,
where $w_{s}$ is a $d$-dimensional Wiener process,
$\sigma$ is a uniformly nondegenerate, bounded,
Borel function with values in the set of symmetric
$d\times d$ matrices, $b$ is a Borel measurable $\bR^{d}$-
valued function given on $\bR^{d+1}=(-\infty,\infty)\times\bR^{d}$
such that
\begin{equation}
                                                    \label{5.10.2}
\int_{\bR}\Big(\int_{\bR^{d}}|b(t,x)|^{p_{0}}\,dx\Big)^{q_{0}/p_{0}}\,dt<\infty 
\end{equation}
if $p_{0}\geq q_{0}$ or
\begin{equation}
                                                    \label{5.10.20}
\int_ {\bR^{d}} \Big(\int_{\bR}|b(t,x)|^{q_{0}}\,dt\Big)^{p_{0}/q_{0}}\,dx<\infty 
\end{equation}
if $ p_{0}\leq q_{0}$. 
Observe that the case $p_{0}=q_{0}=d+1$ is not excluded and in this case
the condition becomes $b\in L_{d+1}(\bR^{d+1})$.

The goal of this article is to study some properties of such
solutions or Markov processes whose trajectories
are solutions of \eqref{11.29.2}. In particular,
in Section \ref{section 10.25.2} for more or less general
processes of diffusion type we derive
several estimates of  Aleksandrov type by using
  Lebesgue spaces with mixed norms
like (in case $t=0$, $x=0$ in \eqref{11.29.2})
\begin{equation}
                                               \label{10.26.2}
E\int_{0}^{\infty}e^{- t}
f(t,x_{t}) \,dt\leq 
N \|  f\|_{L_{p ,q }},
\end{equation}
provided that $d/p+1/q\leq1$.

We also
show that expected time when the process $(t,x_{t})$,
starting at $(0,0)$, exits from $[0,R^{2})\times\{x:|x|<R\}$
is comparable to $R^{2}$. This plays a crucial role
in Section \ref{section 10.26.1} where we show a significant 
improvement of the Aleksandrov estimates in the direction
of lowering the powers of summability of $f$ in \eqref{10.26.2}
to $d_{0}/p+1/q\leq1$ with $d_{0}<d$. Time homogeneous
versions of these estimates are also given.

In the same Section \ref{section 10.25.2} we give some estimates
of the distribution of the exit times from cylinders,
which are heavily used in the sequel. We also prove
that, for any $0\leq s\leq t<\infty$,
$$
E\sup_{r\in[s,t]}|x_{r}-x_{s}|^{ n}
\leq N(|t-s|^{ n/2}+|t-s|^{n}).
$$
 
It is to be said that instead of \eqref{5.10.2}
or \eqref{5.10.20}, which are not invariant
under self-similar transformations, we impose
a slightly stronger assumption on $b$,
that is invariant.

As we mentioned above, in Section \ref{section 10.26.1}
we improve the results of Section \ref{section 10.25.2}
in what concerns the Aleksandrov estimates, which allows us
to prove It\^o's formula for $ W^{1,2}_{p,q}(Q)$-functions
if $d_{0}/p+1/q\leq1$.

In Section \ref{section 10.13.1} we discuss some applications
of our results to the theory of parabolic
equations. It\^o's formula is the main instrument here.
 We prove the qualitative form of the parabolic Aleksandrov
maximum principle for $u\in W^{1,2}_{p,q}$ with 
$d_{0}/p+1/q\leq1$ (Theorem \ref{theorem 10.14.1}).
In the case of bounded $b$ and $p=q=d+1$
in the parabolic case and $p=d$ in the elliptic case
the result of Theorem \ref{theorem 10.14.1} is ``classical''
(about 50 year old). It
was 
generalized by Cabr\'e \cite{Ca_95},
 Escauriaza \cite{Es_93}, and Fok \cite{Fo_95} in the elliptic
case when $p<d$ (close to $d$) again
  when $b$ is   bounded. In \cite{CKS_00}
a parabolic version of these results, extending some earlier
results by Wang, are given
for $L_{p}$-viscosity solutions with
$p<d+1$ (close to $d+1$) when $b$ is   bounded.
However, it is  worth noting that in the elliptic case
it may happen that $b\not\in L_{d}$
and the equation is still solvable (see, for instance,
\cite{KS_19}).
 In our situation
we have some freedom in choosing $p,q$ and $b\in L_{p_{0},q_{0}}$,
but we only treat true solutions.
Theorem \ref{theorem 10.14.1} covers
Theorem 2.4 of \cite{CKS_00} on the account of
having mixed norms
and $b\in L_{p_{0},q_{0}}$.

One more results in this section is
 aimed at applications to the theory
of fully nonlinear parabolic equations 
with lower order coefficients in $L_{p,q}$.  
We prove  a theorem allowing one
to pass to the limit under the sign of fully
nonlinear operator when the arguments (functions)
 converge only weakly and give its application
to linear equations.

It is worth mentioning that there is a vast 
literature about stochastic equations when
\eqref{5.10.1} is replaced with $d/p+2/q\leq1$.
This condition is much stronger than ours.
Still we refer the reader to the recent articles
\cite{Na_18},
\cite{BFGM_19},  \cite{XZ_20}  and the 
references therein
for the discussion of many powerful 
and exciting results obtained
 under this stronger condition.
There are also many papers when this condition is
considerably relaxed on the account of
imposing various regularity conditions
on $\sigma$ and $b$ and/or considering
random initial conditions with bounded density,
 see, for instance,
\cite{ZZ_19}, \cite{Zh_20} and the 
references therein. Restricting
the situation to the one when $\sigma$ and $b$
are independent of time allows one to
relax the above conditions significantly
further, see, for instance, \cite{KS_19}
and the 
references therein.

 Introduce
$$
B_{R} =\{x\in\bR^{d}:| x|<R\}, 	\quad B_{R}(x)=x+B_{R},
\quad C_{T,R}=[0,T)\times B_{R},
$$
$$
C_{T,R}(t,x)=(t,x)+C_{T,R},\quad C_{R}(t,x)=C_{R^{2},R}(t,x),
\quad C_{R} =C_{R}(0,0),
$$
$$
 D_{i}=\frac{\partial}{\partial x^{i}},
\quad D_{ij}=D_{i}D_{j}\quad \partial_{t}=\frac{\partial}{\partial t}.
$$
For $p,q\in[1,\infty]$ and domains $Q\subset \bR^{d+1}$
 we introduce the space $L_{p,q}(Q)$ as the space of Borel
functions on $Q$ such that
$$
\|f\|^{q}_{L_{p,q}(Q)}:=
\int_{\bR}\Big(\int_{\bR^{d}}I_{Q}(t,x)|f(t,x)|^{p}\,dx\Big)^{q/p}\,dt<\infty 
$$
if $p\geq q$ or
$$
\|f\|^{p}_{L_{p,q}(Q)}:=\int_ {\bR^{d}} \Big(\int_{\bR}I_{Q}(t,x)
|f(t,x)|^{q}\,dt\Big)^{p/q}\,dx<\infty 
$$
if $ p\leq q$ with natural interpretation
of these definitions if $p=\infty$ or $q=\infty$.
If $Q=\bR^{d+1}$, we drop $Q$ in the above notation.
Observe that
$p $ is associated with   $x$ and
$q$ with   $t$ and the interior
integral is always elevated to the power $\leq 1$.
 In case $p=q=d+1$ we abbreviate $L_{d+1,d+1}=L_{d+1}$.
For the set of functions on $\bR^{d}$ summable to the $p$th
power we use the notation $L_{p}(\bR^{d})$.

If $\Gamma$ is a measurable subset of $\bR^{d+1}$ we denote by
$|\Gamma|$ its Lebesgue measure. The same notation
is used for measurable  subsets of $\bR^{d}$ with $d$-dimensional
Lebesgue measure. We hope that it will be clear
from the context which Lebesgue measure is used.
If $\Gamma$ is a measurable subset of $\bR^{d+1}$ and
$f$ is a function on $\Gamma$ we denote 
$$
\dashint_{\Gamma}f\,dxdt=\frac{1}{|\Gamma|}
\int_{\Gamma}f\,dxdt.
$$
In case $f$ is a function on a measurable subset $\Gamma$
of $\bR^{d}$ we set
$$
\dashint_{\Gamma}f\,dx =\frac{1}{|\Gamma|}
\int_{\Gamma}f\,dx .
$$

Throughout the article $\bar R$ is 
a fixed constant, $\bar R\in(0,\infty)$.

\mysection{The case of general diffusion type processes
with drift in $L_{p_{0},q_{0}}$}

                                        \label{section 10.25.2}

Let $(\Omega,\cF,P)$ be a complete probability
space, let $\cF_{t}, t\geq0$, be an increasing family of
complete $\sigma$-fields $\cF_{t}\subset\cF$,  
let $w_{t}$ be an $\bR^{d}$-valued Wiener process
relative to $\cF_{t}$. Fix $\delta\in (0,1)$
and denote by $\bS_{\delta}$
the set of $d\times d$ symmetric matrices
whose eigenvalues are between $\delta$ and
$\delta^{-1}$. Assume that we are given
an $\bS_{\delta}$-valued $\cF_{t}$-adapted process
$\sigma_{t}=\sigma_{t}(\omega)$ and an $\bR^{d}$-valued
$\cF_{t}$-adapted process $b_{t}$, such that
$$
\int_{0}^{T}|b_{t}|\,dt<\infty
$$
for any $T\in(0,\infty)$ and $\omega$. Define
$$
x_{t}=\int_{0}^{t}\sigma_{s}\,dw_{s}+\int_{0}^{t}b_{s}\,ds.
$$

\begin{assumption}
                                    \label{assumption 10.20.1}
We are given a function $h\in L_{p_{0},q_{0},\loc}$
such that
$$
|b_{t}|\leq  h(t
,x_{t}).
$$
Furthermore, there exists a  {\em
bounded\/} nondecreasing  function 
$\bar b_{R}$, $R\in(0,\infty)$, such that
 for any $(t,x)\in\bR^{d+1}$ 
    and  $R\in(0,\infty)$ we have
\begin{equation}
                                               \label{8.19.1}
 \|h\|^{q_{0}}_{L_{p_{0},q_{0}}(C_{R}(t,x))} \leq
\bar b_{R}R .
\end{equation} 
\end{assumption}

Observe that if $p_{0}=q_{0}=d+1$ and $h\in L_{d+2}$
(a typical case in the theory of parabolic equations),
then \eqref{8.19.1} is satisfied with $\bar b_{R}=\bar b=
\|h\|^{d+1}_{L_{d+2}}$, because by H\"older's inequality
$$
\|h\|^{d+1}_{L_{d+1}(C_{R}(t,x))}
\leq \bar bR.    
$$
 On the other hand, it may happen that
\eqref{8.19.1} is satisfied with $p_{0}=q_{0}=d+1$ 
but $h\not\in L_{d+2,\loc}$.

\begin{example}
                                          \label{example 8.24.1}
Take  $\alpha\in(0,d),\beta\in(0,1)$ such that
$\alpha+2\beta=d+1$ and consider the function
$g(t,x)=|t|^{-\beta}|x|^{-\alpha}$. Observe that
$$
\int_{C_{R}(t,x)}g(s,y)\,dyds  
=R  \int_{C_{1}(t',x')}g(s,y)\,dyds,
$$
where $t'=t/R^{2}$, $x'=x/R$. Obviously, the last integral
is a bounded function of $(t',x')$. Hence, the function
$h=g^{1/(d+1)}$ satisfies \eqref{8.19.1}
with $p_{0}=q_{0}=d+1 $. As is easy to see
for any $p>d+1$ one can find $\alpha$ and $\beta$
above such that $h\not\in L_{p,\loc}$.
 
\end{example}

Note that if $h$ is bounded and has compact support, \eqref{8.19.1}
is certainly satisfied.
A condition very similar to
\eqref{8.19.1} first appeared in \cite{Na_12}. 

The following is a particular case 
of Theorem 4.5 of \cite{Kr_20_2}. In Theorem \ref{theorem 9.27.1}
$p =\infty$ is allowed (and then $q =1$).

\begin{theorem}
                                        \label{theorem 9.27.1}

Suppose that  Assumption \ref{assumption 10.20.1}
is satisfied and
\begin{equation}
                                                  \label{10.4.3}
p,q\in[1,\infty],\quad\frac{d}{p }+\frac{1}{q }=1 .
\end{equation}
Then for any Borel $f\geq0$ and stopping time $\gamma$
\begin{equation}
                                                   \label{5.6.4}
  E\int_{0}^{\gamma}   
f( t,x_{t})\,dt\leq N(d,p_{0},\delta  )
\Big(A+ \|h\|_{L_{p_{0},q_{0}}}^{2q_{0} }
\Big)^{d/(2p )}
\|f\|_{L_{p ,q }},
\end{equation}
where $A=E\gamma$.

\end{theorem}

Our first goal is to estimate $A$ and eliminate it from 
\eqref{5.6.4}. For $T,R\in(0,\infty)$ introduce
$$
\tau_{T,R}(x)=\inf\{t\geq 0:(t,x+x_{t})\not\in C_{T,R}\},
\quad \tau_{R}(x)=\tau_{R^{2},R}(x),\quad \tau_{R}=\tau_{R}(0) .
$$

\begin{lemma}
                                      \label{lemma 8.16.1}
 
We have
\begin{equation}
                                          \label{8.16.1}
A:=
E \tau_{R} (x)\leq  R^{2},
\end{equation}
 and consequently, assuming \eqref{10.4.3},
 for any Borel nonnegative $f$
\begin{equation}
                                          \label{9.29.2}
 E\int_{0}^{\tau_{R}(x)}  
f( t,x_{t})\,dt\leq 
 N(d,p_{0},  \delta)(1+\bar b_{R}) ^{d/  p  }
 R 
 ^{d/  p  }
\|f\|_{L_{p ,q }}. 
\end{equation}
 
\end{lemma}

Proof. Obviously, $\tau_{R}\leq R^{2}$ 
and \eqref{8.16.1} follows.  After that \eqref{9.29.2}
follows from \eqref{8.19.1} and \eqref{5.6.4}.
The lemma is proved.

Estimate \eqref{8.16.1} says that in the typical case
of nondegenerate diffusion $\tau_{R}$ is of order
not more than $R^{2}$.
A very important fact which is implied by Corollary
\ref{corollary 7.29.1} is that 
$\tau_{R}$ is of order
not less than $R^{2}$. To show this we need
an additional assumption appearing after the following result,
in which
\begin{equation}
                                                  \label{1.6.1}
\tau'_{R}(x)=\inf\{t\geq0:x+x_{t}\not\in B_{R}\},\quad
\gamma_{R}(x)=\inf\{t\geq0:x+x_{t} \in \bar B_{R}\}.
\end{equation}

\begin{theorem}
                                      \label{theorem 8.2.1}
There are
   constants $\bar \xi=\bar \xi(d,\delta)\in (0,1) $ and
 $\bar N=\bar N(d,p_{0},  \delta)$ {\em
continuously\/} depending on $\delta$
such that if, for an $  R\in(0,\infty)$, we have
\begin{equation}
                                     \label{12.18.3}
\bar N \bar b_{  R}\leq 1,
\end{equation}
then  for $|x|\leq R$
\begin{equation}
                                          \label{8.2.2} 
  P( \tau_{R}(x)  =   R^{2} )\leq 1-\bar\xi,
\quad P( \tau_{R}  =   R^{2} )\geq 
 \bar\xi .   
\end{equation}
Moreover for $n=1,2,...$ and $|x|\leq R$
\begin{equation}
                                          \label{1.3.1} 
P(\tau'_{R}(x)\geq nR^{2})
=  P( \tau_{nR^{2},R}(x)  =  n R^{2} )\leq (1-\bar\xi)^{n},   
\end{equation}
so that $E\tau'_{R}(x)\leq N(d,\delta)R^{2}$
and
\begin{equation}
                                          \label{1.3.3} 
I:=E\int_{0}^{\tau'_{R}(x)}h(t,x_{t})\,dt
\leq N(d,p_{0},\delta)
(1+\bar b_{R}) ^{d/  p_{0}  }\bar b^{1/q_{0}}_{R}R.
\end{equation}

Furthermore,  the probability
starting from a point in the closed ball of radius $R/16$
with center in   $\bar B_{ R/2}$
to reach the ball $\bar B_{R/16}$ before exiting
from $B_{R}$ is bigger than $\bar\xi$:
for any $x,y$ with $|y|\leq R/2$ and $|x-y|\leq R/16$
\begin{equation}
                                          \label{1.2.1} 
P(\tau'_{R}(x)>\gamma_{R/16}(x))\geq\bar\xi.   
\end{equation}

\end{theorem}

\begin{remark}
                                       \label{remark 1.23.1}
The last statement of the theorem might look
awkward because it just says that for any $x\in \bar B_{9R/16}$
estimate \eqref{1.2.1} holds. Mentioning $y$ seems superfluous.
The goal of introducing $y$ is that \eqref{1.2.1}
shows that starting from any point in $B_{R/16}(y)$
the process reaches $\bar B_{R/16}$ with positive probability
without exiting from $B_{R}$, thus ``moves in the direction''
of $-y$,
 no matter where in $B_{R/16}(y)$ the starting point is.

\end{remark}

We first prove an auxiliary result, in which
$$
m_{t}=\int_{0}^{t}\sigma_{s}\,dw_{s},\quad a_{t}=(1/2)
\sigma_{t}\sigma_{t}^{*}.
$$

\begin{lemma}
                                       \label{lemma 1.2.1}
(i) There exists $\kappa=\kappa(d)>0$ such that
for
$$
\psi(x,t)=R^{-4}\big(R^{2}-4|x|^{2}\big)^{2}\phi_{t},\quad
\phi_{t}=\exp\int_{0}^{t}
\kappa R^{-2}\,\tr a_{s}\,ds
$$
the process $\psi(m_{t},t)$ is a local submartingale.

(ii) Take a $\zeta\in C^{\infty}_{0}(\bR)$
such that it is even, nonnegative, and decreasing
on $(0,\infty)$.
For  $T\in(0,\infty)$ and $x\in \bR$ and $t\leq T$define
$u(t,x)=E\zeta(x+w^{1}_{T-t}) )$. Also 
take $x\in\bR^{d}$ and set
$$
r_{t}=\frac{(x+m_{t},a_{t}(x+m_{t}))}
{|x+m_{t}|^{2}}\quad (0/0:=1),\quad
\eta_{t}=2\int_{0}^{t}r_{s}\,ds.
$$
Then  the process $u(\eta_{t},
|x+m_{t}|)$
is a supermartingale before $\eta_{t}$ reaches $T$,
in particular, on $[0,\delta^{2} T]$.

(iii) There exists $\alpha=\alpha(d,\delta)>1$
such that for $u(x)=|x|^{-\alpha}$
and any nonzero $x\in\bR^{d}$ the process
$u(|x+m_{t}|)$ is a submartingale
before $x+m_{t}$ hits the origin.

\end{lemma}

Proof. It is easy to see that for a $\kappa=\kappa(d)>0$
we have $\kappa\mu^{2}-16\mu+32d^{-1}(1-\mu)\geq0
$ for all $\mu$, which implies that
for all $\lambda$
\begin{equation}
                              \label{12.26.1}
\kappa(1-4\lambda^{2})^{2}-16(1-4\lambda^{2})
+128 d^{-1} \lambda^{2}\geq0.
\end{equation}

It follows that
$$
R^{4}\phi_{t}^{-1}d\psi(m_{t},t)=\kappa 
\big(R^{2}-4|m_{t}|^{2}\big)^{2}R^{-2}\tr a_{t}\,dt
$$
$$
-8\big(R^{2}-4|m_{t}|^{2}\big)\big(2m_{t}\,dm_{t}
+2\tr a_{t}\,dt\big)
+128(m_{t},a_{t}m_{t})\,dt\geq dM_{t},
$$
where $M_{t}$ is a local martingale. This proves (i).

(ii) Observe that $u$ is smooth, even in $x$, and satisfies
$\partial_{t}u+(1/2)u''=0$.
Furthermore, as is easy to see $u'(t,x)\leq 0$
for $x\geq0$.
It follows by It\^o's formula that 
before $\eta_{t}$ reaches $T$ we have
(dropping obvious values of some arguments)
$$
du(\eta_{t},|x+m_{t}|)= r_{t} (2\partial_{t}u+u'')
\,dt +\frac{u'}{|x+m_{t}|} (\tr a_{t}-
 r_{t} )\,dt
+dM_{t},  
$$
where $M_{t}$ is a stochastic integral. Here the second
term with $dt$ is negative since $u'\leq0$, and this proves 
that $u(\eta_{t},|x+m_{t}|)$ is a local supermartingale.
Since it is nonnegative, it is a supermartingale.

Assertion (iii) is proved by simple application
of It\^o's
formula (see, for instance,
the proof of Lemma 2.2 in \cite{Sa_10}).
The lemma is proved.

{\bf Proof of Theorem \ref{theorem 8.2.1}}. 
 Notice that by 
 \eqref{9.29.2}    
\begin{equation}
                                              \label{1.3.4}
E\int_{0}^{\tau_{R}}|b(t,x_{t})| \,dt
\leq  N(1+\bar b_{R})^{d/  p_{0}  }
 \bar b_{R}^{1/q_{0}}  R .
\end{equation}

Furthermore, observe that for $\gamma$ defined as the minimum of
$R^{2}$ and the first exit time of
$m_{t}$ from $B_{R/2}$ it holds that
$
\phi_{\gamma}\leq e^{\kappa d/\delta}.
$
Hence, by Lemma \ref{lemma 1.2.1} (i)
$$
1=\psi(0,0)\leq E\psi(m_{\gamma},\gamma)
\leq e^{\kappa d/\delta^{2}}
P(\sup_{t\leq R^{2}}|m_{t}|< R/2)
$$
and, since $\tau_{R}\leq R^{2}$,
$$
P(\sup_{t\leq \tau_{R}}|m_{t}|< R/2)
\geq 2\bar \xi(d,\delta)>0.
$$
Also note that
$$
P( \tau_{R}  < R^{2} )\leq
P\big(\int_{0}^{\tau_{R}}|b(\sft_{t},x_{t})| \,dt
\geq R/2\big)+
P(\sup_{t\leq \tau_{R}}|m_{t}|\geq R/2).
$$
Therefore we, get the right estimate in \eqref{8.2.2}  for 
 $2N(1+\bar b_{R})^{d/  p_{0}  } 
\bar b_{R}^{1/q_{0}}\leq 
\bar\xi$.

On the other hand,  take $\zeta$ such that
$\zeta(x)=\eta(x/R)$, where $\eta(x)=1$
 for $|x|\leq 2 $ and take
$T=\delta^{2}R^{2} $, in which case $u(0,x)\leq u(0,0)<1$
and $u(0,0)$ depends only on $\delta$ (and $\eta$).
Also define $\mu$ as the first time $\eta_{t}$
reaches $T$, which is certainly less than or
equal to $R^{2}$. Now observe that
$u(\eta_{\mu},|x+m_{\mu}|)=u(T,|x+m_{\mu}|)=\zeta
(|x+m_{\mu}|)$. It follows that
 $$
 P(\sup_{t\leq R^{2}}|x+m_{t}|< 2R)
\leq P( |x+m_{\mu}|< 2R)
$$
$$
\leq Eu(\eta_{\mu},|x+m_{\mu}|)\leq u(0,x)\leq u(0,0).
$$
Hence,
$$
P( \tau_{R}(x) < R^{2} )\geq
 P\big(\int_{0}^{\tau_{R}}|b(\sft_{t},x_{t})| \,dt
\leq R/2,\sup_{t\leq R^{2}}|x+m_{t}|\geq 2R \big)
$$
$$
\geq 1-P\big(\int_{0}^{\tau_{R}}|b(\sft_{t},x_{t})| \,dt
\geq R/2\big)-P(\sup_{t\leq R^{2}}|x+m_{t}|\leq 2R)
$$
and it is clear how to adjust \eqref{12.18.3}
to get both inequalities in \eqref{8.2.2} with perhaps
$\bar\xi$ different from the above one.
Estimate \eqref{1.3.1} is obtained by iterations.

To prove \eqref{1.3.3} come back to \eqref{1.3.4}
and denote by $J$ its right-hand side.   
Then use the condition version of \eqref{1.3.4}
to see that ($\tau_{0,
 R}:=0$)
$$
I=\sum_{n=1}^{\infty}EI_{\tau_{(n-1)R^{2},R}(x)
>\tau'_{R}(x)}E\Big(
\int_{\tau_{(n-1)R^{2},R}(x)}^{\tau_{nR^{2},R}(x)}
h(t,x_{t})\,dt\mid \cF_{\tau_{(n-1)R^{2},R}(x)}\Big)
$$
$$
\leq J\sum_{n=1}^{\infty}P(\tau_{(n-1)R^{2},R}(x)
=(n-1)R^{2})\leq J\sum_{n=1}^{\infty}
(1-\bar\xi)^{n-1}.
$$
This yields \eqref{1.3.3}.  

To prove \eqref{1.2.1} 
use assertion (iii) of Lemma \ref{lemma 1.2.1}
to conclude that
$$
du(|x+x_{t}|)\geq b^{i}_{t}D_{i}u(|x+x_{t}|)\,dt+dM_{t},
$$
where $M_{t}$ is a local martingale. 
For our $x$, on the time interval, which we denote
 $(0,\nu)$, when $x+x_{t}
\in B_{R}\setminus \bar B_{R/16}$ we have
$|D u(|x+x_{t}|)t\leq N(d,\alpha)R^{-\alpha-1}$
Furthermore, at starting point $u(x)\geq (9R/16)^{-\alpha}$.
Consequently and by \eqref{1.3.3}
$$
(9R/16)^{-\alpha}\leq NR^{-\alpha-1}
E\int_{0}^{\tau'_{R}(x)}h(t,x_{t})\,dt
+P\big(\nu=\tau'_{R}(x)\big)R^{-\alpha}
$$  
$$
+
P\big(\nu=\gamma_{R/16}(x)\big)(R/16)^{-\alpha},
$$
$$
(16/9)^{\alpha}\leq N_{1}
(1+\bar b_{R}) ^{d/  p_{0}  }\bar b^{1/q_{0}}_{R}
+1
$$
$$
-P\big(\tau'_{R}(x)>\gamma_{R/16}(x)\big)+
16^{\alpha}P\big(\tau'_{R}(x)>\gamma_{R/16}(x)\big).
$$
It follows easily that \eqref{1.2.1} holds
with $\bar\xi$ perhaps different from the above ones,
once a relation like \eqref{12.18.3} holds.
The continuity of $\bar N$ in \eqref{12.18.3}
and of $\xi(d,\delta)$
with respect to $\delta$ is established by inspecting
the above proof.
The theorem is proved.

\begin{assumption}
                                   \label{assumption 12.18.2}
There exists $\uR\in(0,\infty)$
such that
\begin{equation}
                                     \label{1.23.2}
\bar N(d,p_{0},\delta) \bar b_{ \uR}< 1.
\end{equation}
\end{assumption}

This assumption as well as Assumption  
\ref{assumption 10.20.1} 
is supposed to hold throughout the article.
 Set
$$
\ulambda=\uR^{-2}.
$$

\begin{corollary}
                                          \label{corollary 8.18.1}
For  $\mu \in[0,1]$ and   $R\leq\uR$ we have
\begin{equation}
                                          \label{8.18.3}
Ee^{-\mu R^{-2} \tau_{R}} \leq e^{- 
\mu
\bar \xi/2}.
\end{equation}
\end{corollary} 

 Indeed,
the derivative with respect to $\mu$ of the left-hand
side of \eqref{8.18.3} is
$$
-R^{-2} E \tau_{R} e^{-R^{-2}\mu \tau_{R}} \leq -
 e^{-\mu }R^{-2}
P( \tau_{R}=R^{2})\leq-e^{-\mu}\bar \xi,
$$
where the last inequality follows from \eqref{8.2.2}.
By integrating we find
$$
Ee^{-\mu R^{-2} \tau_{R}} -1\leq  
(e^{-\mu   }-1)
\bar \xi,
$$
which after using
$$
e^{-\mu  }-1\leq- \mu /2,\quad
1-\mu
\bar \xi/2\leq e^{-\mu 
\bar \xi/2}
$$
leads to \eqref{8.18.3}.

\begin{theorem}
                                     \label{theorem 8.20.1}
For any $\lambda,R>0$ we have
\begin{equation}
                                          \label{8.20.1}
Ee^{-\lambda  \tau_{R} }\leq
e^{\bar\xi/2}e^{- \sqrt{\nelambda}  R
\bar\xi/2}=
\begin{cases}
 e^{\bar\xi/2}e^{-\sqrt\lambda R\bar\xi/2}
\quad\text{if}\quad \lambda
 \geq \ulambda\\
e^{\bar\xi/2}e^{-\lambda R\uR \bar\xi/2}\quad\text{if}
\quad \lambda
 \leq \ulambda,
\end{cases}
\end{equation}
where 
$$
\nelambda= \lambda\min(1, \lambda/\ulambda ).
$$ 
In particular, for  
  any   $R>0$ and $t\leq R\uR
\bar \xi/4 $ we have
\begin{equation}
                                             \label{10.2.2}
P(  \tau_{R}  \leq  t )\leq 
 e^{\bar\xi/2}\exp\Big(-\frac{{\bar \xi}^{2}R^{2}}{16 t}\Big).
\end{equation}

\end{theorem}

Proof. 
Take an integer $n\geq 1$, introduce $\tau^{k}$, $k=1,...,n$,
as the first exit time of $(\sft_{t},x_{t})$
from $C_{R/n}( \tau^{k-1} ,x_{\tau^{k-1}})$ after $\tau^{k-1}$
($\tau^{0}:=0$). 
If
$$
\lambda\leq n^{2}/R^{2},\quad R/n\leq \uR,
\quad\text{that is}\quad
 n\geq R\sqrt\lambda\max(1,/(\sqrt\lambda \uR)),
$$
then by \eqref{8.18.3} with $\mu=(R/n)^{2}
\lambda$ we have
$$
E\Big(e^{-\lambda( \tau^{k} - \tau^{k-1} )}\mid \cF_{\tau^{k-1}}\Big)
\leq e^{- (R/n)^{2}
\lambda \bar \xi/2}.
$$
Hence,
\begin{equation}
                                           \label{12.27.1}
Ee^{-\lambda \tau_{R} }\leq E\prod_{k=1}^{n}
e^{-\lambda( \tau^{k} - \tau^{k-1}  )}\leq  
 e^{- R^{2}n^{-1}
\lambda \bar \xi/2}.
\end{equation}
If $\lambda \uR^{2}\geq1$, we take
$n=\lceil R\sqrt\lambda\rceil$ and use that
$R^{2}n^{-1}
\lambda\geq R\sqrt\lambda-1$.
If $\lambda \uR^{2}\leq1$, we take $n=\lceil R/\uR\rceil$
and use that
$R^{2}n^{-1}
\lambda\geq R\uR \lambda-1$.
This proves \eqref{8.20.1}.

To prove \eqref{10.2.2} 
observe that if $\lambda  \geq \ulambda$
$$
P( \tau_{R}  \leq  t  )=P\big(
 \exp(-\lambda 
 \tau_{R} ) \geq \exp(-\lambda t)\big)
\leq  
 \exp(\bar \xi/2+ \lambda t-\sqrt{\lambda}R \bar \xi/2).
$$
For $\sqrt{\lambda}=R\bar \xi/(4t)$
 we get \eqref{10.2.2} provided
$R\bar \xi/(4t)\geq \uR^{-1}$.
The theorem is proved.

Recall that $\bar R$ is fixed throughout the article.

\begin{corollary}
                                       \label{corollary 7.29.1}
Let   $  \Lambda\in(0,\infty)$.
Then there is a  constant   
  $N =N (\uR,\bar R,\Lambda,\bar\xi)$
such that for any $R\in(0,\bar R]$, $\lambda\in[0,\Lambda]$ 
\begin{equation}
                                                   \label{8.21.1}
N E \tau_{R} \geq R^{2},\quad N E\int_{0}^{\tau_{R}}
e^{-\lambda t} \,dt  \geq R^{2} .
\end{equation}

\end{corollary}

Indeed, for any $\nu\leq  \uR\bar \xi/(4\bar R)$
and $R\in(0,\bar R]$
we have $\nu R^{2}\leq R\uR\bar\xi/4$ so that
$$
E \tau_{R} \geq  \nu
R^{2}P( \tau_{R} >\nu R^{2})\geq 
\nu R^{2}
 \Big(1-e^{\bar\xi/2} \exp
\Big(-\frac{\bar\xi^{2} }{16 \nu}\Big)\Big),
$$
$$
E\int_{0}^{\tau_{R}}
e^{-\lambda t} \,dt=\lambda^{-1}
E(1-e^{-\lambda \tau_{R }})
\geq\lambda^{-1}
 EI_{ \tau_{R}  >\nu R^{2}} (1-e^{-\lambda\nu R^{2}})
$$
$$
=\lambda^{-1} P( \tau_{R} >\nu R^{2})(1-e^{-
\lambda\nu R^{2}})
$$
$$
\geq \lambda^{-1}\Big(1-e^{\bar\xi/2} \exp
\Big(-\frac{\bar\xi^{2} }{16\nu }\Big)\Big)
(1-e^{-\lambda\nu R^{2}}),
$$
which yields \eqref{8.21.1} for an appropriate 
small 
$\nu =\nu(\uR,\bar R,\Lambda,\bar\xi )>0$.

\begin{corollary}
                                        \label{corollary 10.26.1}
For any $n>0$ and
  $0\leq s\leq t $ we have
\begin{equation}
                                            \label{10.28.2}
E\sup_{r\in[s,t]}|x_{r}-x_{s}|^{ n}
\leq N(|t-s|^{ n/2}+|t-s|^{ n}),
\end{equation}
where $N=N(n,\uR,\bar\xi)$.
\end{corollary}

Indeed, clearly we may assume that $s=0$. Then
for $\nu_{0}=4(\uR\bar\xi)^{-1}$ and $\mu\geq
t\nu_{0}$ we have $t\leq \mu\uR\bar\xi/4$
and 
$$
P(\sup_{r\leq t }|x_{r} |\geq \mu)
\leq P(\tau_{\mu}\leq t)
\leq e^{\bar\xi/2}\exp\Big(-\frac{  \mu^{2}\bar\xi^{2}}
{16 t}\Big).
$$
Consequently,
$$
E\sup_{r\leq t]}|x_{r} |^{ n}
=n\int_{0}^{\infty}\mu^{n-1}
P(\sup_{r\leq t }|x_{r} |\geq \mu)\,d\mu
\leq n\int_{0}^{t\nu_{0}}\mu^{n-1}\,d\mu
$$
$$
+ne^{\bar\xi/2}\int_{0}^{\infty}\mu^{n-1}
\exp\Big(-\frac{  \mu^{2}\bar\xi^{2}}
{16 t}\Big)
\,d\mu,
$$
and   the result follows. 

A few more general results are related to going through
a long ``sausage".
\begin{theorem}
                                        \label{theorem 1.24.1}
Let $R\in(0,\uR]$, $x,y\in\bR^{d}$ and $16|x-y|\geq 3R$.  
For $r>0$ denote by $S_{r}(x,y)$   the open convex hull
of $B_{r}(x)\cup B_{r}(y)$. Then there exist
$T_{0},T_{1}$, depending only on $\bar\xi$,
such that $0<T_{0}<T_{1}<\infty$ and the probability $\pi$
that $x+x_{t}$ will reach $\bar B_{R/16}(y)$ before exiting
from $S_{R}(x,y)$ and this will happen
on the time interval $[nT_{0}R^{2},nT_{1}R^{2}]$
is greater than $\pi_{0}^{n}$, where
$$
n= \Big\lfloor \frac{16|x-y|+R}{4R}\Big\rfloor 
$$  
and $\pi_{0}=\bar\xi/3$.

\end{theorem}  

Proof. We may assume that $y=0$.
Introduce $\tau(x)$ as the first time $x+x_{t}$
reaches $\bar B_{R/16}$ and $\gamma(x)$ as the first time
it exits from $S_{R}(x,0)$. Owing to
$16|x |\geq 3R$, we have $n\geq1$ and we are going to use the induction
on $n$ with the induction hypothesis that
$$
\Big\lfloor \frac{16|x |+R}{4R}\Big\rfloor=n
\Longrightarrow 
P(\gamma(x)>\tau(x)\in[nT_{0}R^{2},n
T_{1}R^{2}])\geq \pi^{n}_{0}.
$$

   If $n=1$, $3R/16\leq |x|< 7R/16$ and by Theorem
\ref{theorem 8.2.1} we have $P(\tau'_{R}(x)>\tau(x))\geq\bar\xi$.
Furthermore, in light of  Theorem
\ref{theorem 8.2.1}, there is $T_{1}=T_{1}(\bar\xi)$
such that $P(\tau'_{R}(x)>T_{1}R^{2})\leq \bar\xi/3$.
Using \eqref{10.2.2} we also see that there is
$T_{0}=T_{0}(\bar\xi)<T_{1}$ such that
$P(\tau(x)\leq T_{0}R^{2})\leq \bar\xi/3$.
Hence, $P(\gamma(x)>\tau(x)\in[T_{0}R^{2},
T_{1}R^{2}])\geq\bar\xi/3=\pi_{0}$.
This justifies the start of the induction.

Assuming that our hypothesis is true for some $n\geq 1$
 suppose that
$(n+2)R/4>|x|+R/16\geq  (n+1)R/4$. In that case, let $z=nRx/(4|x|)$,
$\tau_{z}$ be the first time $x+x_{t}$ reaches $\bar B_{R/16}(z)$,
and let $\gamma_{z}$ be the first time it exits
from $S_{R}(x,z)$. As is easy to see
$$
P(\gamma(x)>\tau(x)\in[(n+1)T_{0}R^{2},(n+1)
T_{1}R^{2}])
$$
$$
\geq P(\gamma_{z}>\tau_{z}\in[T_{0}R^{2},T_{1}R^{2}],
\gamma(x_{\tau_{z}})>\tau(x_{\tau_{z}})\in[nT_{0}R^{2},n
T_{1}R^{2}])
$$
$$
=EI_{\gamma_{z}>\tau_{z}\in[T_{0}R^{2},T_{1}R^{2}]}
P\Big(\gamma(x_{\tau_{z}})>\tau(x_{\tau_{z}})\in[nT_{0}R^{2},n
T_{1}R^{2}]\mid \cF_{\tau_{z}}\Big).
$$
Observe that on the set $\tau_{z}<\infty$ we have
$nR/4\leq |x_{\tau_{z}}|+R/16<(n+1)R/4$, so that by the conditional version
of our induction hypothesis the conditional probability
above is greater than $\pi_{0}^{n}$. Then just by shifting the origin
to $z$ and using the first part of the proof
we obtain our result for $n+1$ in place of $n$.
The theorem is proved.

\begin{remark}
                                      \label{remark 1.25.1}
Observe that, for any fixed $x,y$, the interval
$[nT_{0}R^{2},nT_{1}R^{2}]$ is as close to zero
as we wish if we choose $R$ small enough.
Then, of course, the corresponding probability will be
quite small but $>0$.

\end{remark} 
\begin{corollary}
                                     \label{corollary 1.25.1}
Let $R\leq\uR$, $\kappa\in[0,1)$, and $|x|\leq\kappa R$.
Then for any $T>0$
\begin{equation}
                                        \label{1.25.2}
NP(\tau'_{R}(x)> T)\geq e^{-\nu T/[(1-\kappa)R]^{2}},
\end{equation}
where $N$ and $\nu>0$ depend only on $\bar\xi$.
\end{corollary}

Indeed, passing from $B_{R}$ to $B_{(1-\kappa)R}(x)$
shows that we may assume that $x=0$ and $\kappa=0$.
In that case, consider meandering of $x_{t}$
between $\bar B_{R/16}$ and $\partial B_{R/16}(y)$ where
$|y|=R/4$ without exiting from $B_{R}$. As is easy to deduce
from Theorem \ref{theorem 1.24.1},
given that the $n$th loop happened, with probability
$\pi^{4}_{0}$ the next loop will occur and take at least
$4R^{2}T_{0}$ of time. Thus  the $n$th loop
will happen and will take   at least $4nR^{2}T_{0}$ of time
with probability at least $\pi_{0}^{4n}$.
It follows that, for any $n$,
$$
P(\tau'_{R} \geq 4nR^{2}T_{0}))\geq \pi_{0}^{4n},
$$
and this yields \eqref{1.25.2} for $x=0$ and $\kappa=0$.

The following   complements Corollary \ref{corollary 1.25.1}.
\begin{corollary}
                                  \label{corollary 2.3.1}
Let   $R\in(0,\bar R]$. Then there
exists a constant $N$, depending only on
$\bar\xi,\bar R,\uR$, such that, for any 
$T>0$,
$$
P(\tau'_{R}>T)\leq Ne^{-T/(NR^{2})}.
$$
\end{corollary}

Indeed, if $R\leq \uR$, the result follows from Theorem
\ref{theorem 8.2.1}. For $R\geq \uR$, take a point  $y$
such that $|y|=\uR+\bar R$, for any $x$
define $\gamma(x)$ as the first time $x+x_{t}$
hits $\bar B_{\uR/16}(y)$, and set
$$
n_{0}=\Big\lfloor \frac{16(\uR+2\bar R)+\uR}{4\uR}\Big\rfloor.
$$
It follows from Theorem \ref{theorem 1.24.1}
that for any $x\in B_{  R}$
$$
P(\tau'_{R}(x)\leq n_{0}T_{1}\uR^{2})
\geq P(\gamma(x)\leq n_{0}T_{1}\uR^{2})\geq \pi_{0}^{n_{0}}.
$$
Hence 
$$
P(\tau'_{R}(x)> n_{0}T_{1}R^{2})\leq
  1-\pi_{0}^{n_{0}}
$$
and the result follows from Khasminski's lemma.

\mysection{Mixed norm estimates of potentials
of stochastic processes}

Here we are moving toward estimating the resolvents
of Markov diffusion processes in $L_{p,q}$.
\begin{lemma}
                                            \label{lemma 8.22.1}
Assume \eqref{10.4.3}.
Then there is a constant $N$,
depending only on
$\delta,d,  p_{0} $, and $\bar b_{\infty}$, such that
for any   $t_{0}\geq0$, $x_{0}\in \bR^{d}$, $\lambda>0$,
and Borel nonnegative $f$ vanishing outside 
$C_{\nelambda^{-1/2}}(t_{0},x_{0})$
we have 
\begin{equation}
                                             \label{8.22.10}
 E\int_{0}^{\infty}e^{- \lambda t}
f(t,x_{t})  \,dt\leq 
N\nelambda^{-d/(2p)} 
\Phi_{\lambda}(t_{0},x_{0})\|  f\|_{L_{p ,q }},
\end{equation}
where   $\Phi_{\lambda}
(t,x)=e^{-\sqrt{\nelambda}  (\sqrt t+|x|)\bar\xi/4}$.
\end{lemma}

Proof.  
Fix $\rho=N(\bar\xi)\nelambda^{-1/2} >0$ such that
the right-hand side of \eqref{8.20.1} equals $1/2$ when $R=\rho$.
  Then
introduce
  $\tau^{0}$ as the first time $(t,x_{t})$
hits $\bar C_{\nelambda^{-1/2}}(t_{0},x_{0})$ and set $\gamma^{0}$
as the first time after $\tau^{0}$ the process
$(t,x_{t})$ exits from $C_{\nelambda^{-1},
\nelambda^{-1/2} +\rho }(t_{0},x_{0})$.
We define recursively $\tau^{k}$, $k=1,2,...$, as the first
time after $\gamma^{k-1}$ the process
$(t,x_{t})$
hits $\bar C_{\nelambda^{-1/2}}(t_{0},x_{0})$ and $\gamma^{k}$ as 
the first time after $\tau^{k}$ the process
$(t,x_{t})$ exits from 
$C_{\nelambda^{-1},\nelambda^{-1/2} +\rho }(t_{0},x_{0})$.

 These stopping times are either
infinite or lie between $t_{0}$ and $t_{0}+\nelambda^{-1}$.
Therefore,  the left-hand side of
\eqref{8.22.10} equals
\begin{equation}
                                                \label{8.22.1}
 E
\sum_{k=0}^{\infty} e^{- \lambda \tau^{k} }I_{k},
\end{equation}
where
$$
I_{k}=I_{\tau^{k}>t_{0} }
E\Big(\int_{\tau^{k}\wedge(t_{0}+\nelambda^{-1})}^{ \gamma^{k}
\wedge(t_{0}+\nelambda^{-1})}
e^{-\nelambda( t-  \tau^{k })}
f(t,x_{t}) \,dt \mid \cF_{\tau^{k}}\Big).
$$
Here  on the set where $\tau^{k}>t_{0} $
$$
\int_{\tau^{k}\wedge(t_{0}+\nelambda^{-1})}^{ \gamma^{k}
\wedge(t_{0}+\nelambda^{-1})}
 \,dt
=  \gamma^{k}\wedge(t_{0}+\nelambda^{-1}) 
- \tau^{k}  \leq 
 \nelambda^{-1}.  
$$
Using this after estimating the norm of $h$
in $C_{\nelambda^{-1},\nelambda^{-1/2} +\rho }(t_{0},x_{0})$
we infer    from \eqref{5.6.4}
that
$I_{k}\leq N\nelambda^{-d/(2p)}\|f\|_{L_{  p , q }}$,
where $N=N(d,p_{0},\delta,\bar b_{\infty})$.   

Next,   observe that, if $\sqrt{t_{0}}>|x_{0}| $, then  $\tau^{0}$
is bigger than the first exit time of $(t,x_{t})$
from $C_{\sqrt{t_{0}}}$, and
by Theorem \ref{theorem 8.20.1} 
$$
Ee^{- 
\lambda \tau^{0} }\leq Ne^{-\sqrt{\nelambda} 
\sqrt{t_{0}}\bar\xi/2}.
$$
In case $\sqrt{t_{0}}\leq |x_{0}| $ and $|x_{0}|>\nelambda^{-1/2}$
our $\tau^{0}$
is bigger than the first exit time of $(t,x_{t})$
from $C_{|x_{0}|-\nelambda^{-1/2}}$, and
$$
Ee^{-\lambda  \tau^{0} }
\leq Ne^{-\sqrt{\nelambda}(|x_{0}|-\nelambda^{-1/2})\bar\xi/2}.
$$
The last estimate (with   $N=1$) also holds if 
$|x_{0}|\leq\nelambda^{-1/2}$, so that
in case $\sqrt{t_{0}}\leq |x_{0}| $
$$
Ee^{-\lambda  \tau^{0}} 
\leq Ne^{-\sqrt{\nelambda} |x_{0}|\bar\xi/2}
$$
and we conclude that in all cases
$$
Ee^{- \lambda \tau^{0}} \leq 
Ne^{-\sqrt{\nelambda} (\sqrt{t_{0}}+|x_{0}|)\bar\xi/4}.
$$

Furthermore,  by the choice of $\rho$ and
 Theorem \ref{theorem 8.20.1}
$$
E\Big(e^{-\lambda ( \gamma^{k}
- \tau^{k} )}\mid \cF_{\tau^{k}}\Big)\leq \frac{1}{2},
$$
$$
Ee^{-\lambda  \tau^{k }}=
Ee^{-\lambda  \gamma^{k-1} }E
\Big(e^{-\lambda ( \tau^{k} 
- \gamma^{k-1} )} \mid \cF_{\gamma^{k-1}}\Big)
\leq \frac{1}{2}Ee^{- \lambda \gamma^{k-1} },
$$
so that
$$
Ee^{-\lambda  \tau^{k} }\leq\frac{1}{4}
Ee^{- \lambda  \tau^{k-1}} ,\quad
Ee^{- \lambda  \tau^{k}} \leq 4^{-k}
Ee^{- \lambda \tau^{0} }.
$$

Recalling \eqref{8.22.1} we see that the left-hand side 
of \eqref{8.22.10}
is dominated by
$N\Phi_{\lambda}(t_{0},x_{0})\|f\|_{L_{p ,q }}$  
and the lemma is proved.

 The following theorem shows that the time spent
by $(t,x_{t})$ in cylinders $C_{1}(0,x)$
decays very fast as $|x|\to\infty$.
\begin{theorem}
                                            \label{theorem 8.22.1}
Suppose that    
\begin{equation}
                                              \label{10.21.3}
p ,q \in[1,\infty],\quad
\nu:=1-\frac{d}{p }-\frac{1}{q }\geq0 .
\end{equation}
Then
 there is a constant $N=N(\delta,d, p,q, p_{0},\bar b_{\infty})$
 such that
for any $\lambda>0$ and   Borel nonnegative $f$  
we have 
\begin{equation}
                                                   \label{8.22.4}
I:= E\int_{0}^{\infty}e^{- \lambda t}
f(t,x_{t}) \,dt\leq 
N\lambda^{-\nu}\nelambda^{-d/(2p )} \|\Psi_{\lambda}^{1-\nu} f\|_{L_{p ,q }},
\end{equation}
where $\Psi _{\lambda}(t,x)=\exp(- \sqrt{\nelambda} 
(|x|+ \sqrt t)\bar\xi/16)$.
\end{theorem}

Proof. First assume that $\nu=0$.
Take a nonnegative $\zeta\in C^{\infty}_{0}
(\bR^{d+1})$ of the type
$\nelambda^{ (d+2)/2 }\eta(\nelambda t,\sqrt{\nelambda} x)$
 with support in $C_{\nelambda^{-1/2}}$ and
 unit integral and for $(t,x),(s,y)\in\bR^{d+1}$
set
$$
f_{s,y}(t,x)=f(t,x)\zeta(t-s,x-y).
$$
Clearly, due to Lemma \ref{lemma 8.22.1},
$$
I
=\int_{0}^{\infty}\int_{\bR^{d }}E\int_{0}^{\infty}
e^{-\lambda t}
f_{s,y}(t,x_{t}) \,dt\,dyds
$$
$$
\leq N\nelambda^{-d/(2p)} \int_{0}^{\infty}\int_{\bR^{d }}
\Phi_{\lambda}(s,y)\|  f_{s,y}\|_{L_{p ,q }}\,dyds.
$$
{\em Case $p\geq q$}. Then $q<\infty$ and
we introduce
$$
M_{1}^{1/q-1}= 
\int_{0}^{\infty}\int_{\bR^{d }}\Phi^{q/(2q-2)}_{\lambda}(s,y)
\,dyds,
$$
$$
M_{2}^{q/p-1}=
\int_{0}^{\infty}\int_{\bR^{d }}\Phi^{pq/(4p-4q)}_{\lambda}(s,y)
\,dyds,\quad p\ne q,\quad M_{2}=1,\quad p=q.
$$
It follows by H\"older's 
inequality  that
$$
\nelambda^{ d/(2p)} I
\leq NM_{1}\Big(
\int_{0}^{\infty}\int_{\bR^{d }}
\Phi_{\lambda}^{q /2}(s,y)\int_{0}^{\infty}\Big(\int_{\bR^{d}}
f_{s,y}^{p }(t,x)\,dx\Big)^{q /p }dt
\,dyds\Big)^{1/q }
$$
$$
=NM_{1}\Big(
\int_{0}^{\infty}dt \Big(\int_{0}^{\infty}
\int_{\bR^{d }}\Phi_{\lambda}^{q /2}(s,y)
 \Big(\int_{\bR^{d}}
f_{s,y}^{p }(t,x)\,dx\Big)^{q /p } 
\,dyds\Big)\Big)^{1/q } 
$$

$$
\leq NM_{1}M^{1/q}_{2}\Big(
\int_{0}^{\infty}dt\Big(\int_{0}^{\infty}\int_{\bR^{d }}\int_{\bR^{d }}
\Phi_{\lambda}^{p /4}(s,y)f_{s,y}^{p }(t,x)\,dydsdx 
\Big)^{q /p }\Big)^{1/q }.
$$

We replace $\Phi_{\lambda}^{p /4}(s,y)$ by $\Phi_{\lambda}^{p /4}(t,x)$
taking into account that these values are comparable
as long as $\zeta(t-s,x-y)\ne0$. After that
integrating over $dyds$ and computing $M_{1},M_{2}$
lead immediately to
\eqref{8.22.4}. 

{\em Case $p< q$}. 
It follows by H\"older's inequality that  
$$
\nelambda^{ d/(2p)}I\leq NM_{3}\Big(
\int_{0}^{\infty}\int_{\bR^{d }}
\Phi_{\lambda}^{p /2}(s,y)\int_{\bR^{d}}\Big(\int_{0}^{\infty}
f_{s,y}^{q}(t,x)\,dt\Big)^{p/q  }dx
\,dyds\Big)^{1/p }
$$
$$
\leq NM_{3}M_{4}\Big(\int_{\bR^{d}}dx\Big(
\int_{0}^{\infty}\int_{\bR^{d }}\int_{0}^{\infty}
\Phi_{\lambda}^{q /4}(s,y)f_{s,y}^{q}(t,x)\,dtdyds\Big)^{p/q}
\Big)^{1/p},
$$
where
$$
M_{3}^{1/p-1}= 
\int_{0}^{\infty}\int_{\bR^{d }}\Phi^{p/(2p-2)}_{\lambda}(s,y)
\,dyds,
$$
$$
M_{4}^{p/q-1}=
\int_{0}^{\infty}\int_{\bR^{d }}\Phi^{pq/(4q-4p)}_{\lambda}(s,y)
\,dyds .
$$
This leads   to
\eqref{8.22.4} as above.
The theorem is proved if $\nu=0$.

If $\nu\in(0,1)$,
by H\"older's inequality the left-hand side of
\eqref{8.22.4} is dominated by $I_{1}I_{2}$, where
$$
I_{1}^{1/\nu}
= E\int_{0}^{\infty}e^{-\lambda t} \,dt
= 1/\lambda,
$$
$$
I_{2}^{1/(1-\nu)}=
E\int_{0}^{\infty}e^{-\lambda t}
 f^{1/(1-\nu)}(t,x_{t})\,dt.
$$
Here
$$
\frac{d}{p-p\nu}+\frac{1}{q-q\nu}=1
$$
so that by the case that $\nu=0$
$$
I_{2}^{1/(1-\nu)}\leq N\nelambda^{-d/(2p-2p\nu)}\|\Psi_{\lambda}
f^{1/(1-\nu)}\|_{L_{p-p\nu,q-q\nu}}
$$
$$
=N\nelambda^{-d/(2p-2p\nu)}\|\Psi_{\lambda}^{1-\nu}f\|_{L_{p,q}}
^{1/(1-\nu)}.
$$
This leads   to
\eqref{8.22.4} again. Finally, if $\nu=1$ so that
$p=q=\infty$, the left-hand side of \eqref{8.22.4}
is obviously dominated by $\lambda^{-1}\sup f$, so that
\eqref{8.22.4} holds with $N=1$. The theorem is
proved.

By taking $q=\infty$ and $f(t,x)=f(x)$ we come to the following,
which extends Corollary 2.5 of \cite{Kr_19_1} to the case
of time dependent drift $b\in L_{p_{0},q_{0},\loc}$.
It is further generalized  
by relaxing the restriction on $p$
in Theorem \ref{theorem 8.30.1}.
\begin{corollary}
                                      \label{corollary 10.4.1}
Let $p \in [d,\infty]$.
Then   
for any $\lambda>0$ and   Borel nonnegative $f(x)$  
we have 
\begin{equation}
                                                   \label{10.4.2}
 E\int_{0}^{\infty}e^{- \lambda t}
f( x_{t}) \,dt\leq 
N\lambda^{-1+d/p}\nelambda^{-  d/(2p )}\|\Psi_{\lambda}^{d/p} f\|_{L_{p}(\bR^{d})},
\end{equation}
where $\Psi _{\lambda}( x)=\exp(- \sqrt{\nelambda}  |x|\bar\xi/16)$
and $N=N(\delta,d, p, p_{0},\bar b_{\infty})$.
\end{corollary}  
 
Next results are dealing with the exit times of the process
$x_{t}$ rather than $(t,x_{t})$. We will need them
while showing an improved integrability of Green's functions.

Estimate \eqref{9.29.5} below in case $b$ is bounded
was the starting point for the theory
of {\em time homogeneous\/} controlled diffusion processes
about fifty years ago.
 \begin{lemma}
                                      \label{lemma 9.29.1}
Let $p \in [d,\infty]$. Then for any Borel 
nonnegative $f(x)$, $R\leq \bar R$, 
and $x\in\bR^{d}$  
\begin{equation}
                                          \label{9.29.5}
 E\int_{0}^{\tau'_{R}(x)}   
f( x_{t})\,dt\leq N(\delta,d,\bar b_{\bar R},\bar R,\uR , p_{0})
 R 
 ^{2-d/  p  } 
\|f\|_{L_{p  }(\bR^{d})}.
\end{equation}

\end{lemma}

Proof. Define $q$ from \eqref{10.4.3} and
observe that for $k=1,2,...$,
$\Delta_{k} =[(k-1) R^{2},k R^{2})$ and
$$
f_{k}(t,x):=I_{\Delta_{k}}(t)f(x),
$$ 
according to
 \eqref{5.6.4},   on the set where $ \tau'_{R} (x)
\geq  (k-1)R^{2} $ we have
$$
E\Big(\int_{ (k-1)R^{2} }^{(kR^{2})
\wedge \tau'_{R}(x)} f(x_{t})\,dt
\mid \cF_{ (k-1)R^{2} }\Big)
$$
$$
=E\Big(\int_{ (k-1)R^{2}}^{(kR^{2})
\wedge \tau'_{R}(x)} f_{k}(t,x_{t})\,dt  
\mid \cF_{(k-1)R^{2}}\Big)
$$
$$
\leq N\Big(R^{2}+\|h\|_{L_{p_{0},q_{0}}(\Delta_{k}\times B_{R})}^{2q_{0} }
\Big)^{d/(2p )}R^{2/q}
\|f\|_{L_{p  }(\bR^{d})}
$$
$$
\leq N(1+\bar b^{2}_{\bar R})^{d/(2p)}
R^{2-d/p}\|f\|_{L_{p  }(\bR^{d})}.
$$
It follows that
$$
E \int_{0}^{ \tau'_{R}(x)} f(x_{t})\,dt
=\sum_{k=1}^{\infty}
E I_{\tau'_{R}(x)\geq (k-1)R^{2}}
\int_{(k-1)R^{2} }^{(k R^{2})
\wedge \tau'_{R}(x)} f(x_{t})\,dt
$$
$$
\leq NR^{2-d/p}\|f\|_{L_{p  }(\bR^{d})}\sum_{k=1}^{\infty}
P(\tau'_{R}(x)\geq (k-1)R^{2}).
$$
By Corollary \ref{corollary 2.3.1}
 each of the probabilities in the last sum
is less than $Ne^{-k/N}$ and this proves the lemma.

 \mysection{Green's functions}  
                                         \label{section 10.26.1}

Here is a straightforward consequence
of \eqref{8.22.4}.

\begin{theorem}
                                       \label{theorem 9.3.1}
 Assume \eqref{10.21.3} and take $\lambda>0$. Then
there exists  a constant  
$N=N(\delta,d, p,q, p_{0},\bar b_{\infty})$  and
 a nonnegative Borel function $G_{\lambda}(t,x)$
(Green's function of $(\cdot,x_{\cdot})$)
on $\bR^{d+1}$
 such that $G_{\lambda}(t,x)=0$ for $t\leq0$
and for any Borel nonnegative $f$ given on $\bR^{d+1}$
we have
$$
 E\int_{0}^{\infty}e^{-\lambda t}
f(t,x_{t}) \,d t =\int_{\bR^{d+1}}f(t,x)G_{\lambda}(t,x)\,dxdt,
$$
\begin{equation}
                                                   \label{9.3.5}
 \|\Psi^{\nu-1}_{\lambda} G_{\lambda}\|_{L'_{p ,q } } 
\le N\lambda^{-\nu}\nelambda^{ -d/(2p )},
\end{equation}
where we use the notation
$$
\|u\|_{L'_{p,q} }=\Big(\int_{\bR}\Big(
\int_{\bR^{d}}|u(t.x)|^{p '}
\,dx\Big)^{q '/p '}\,dt\Big)^{1/q '}
\quad\text{if}\quad p\geq q,
$$
$$
\|u\|_{L'_{p,q} }=\Big(\int_{\bR^{d}}\Big(
\int_{\bR }|u(t.x)|^{q '}
\,dt\Big)^{p '/q '}\,dx\Big)^{1/p'}
\quad\text{if}\quad p< q,
$$
and $p '=p /(p -1),q '=q /(q -1)$.
 
\end{theorem}

The highest power of pure ($p =q =d+1$) 
summability of $G_{\lambda}$
 guaranteed by this theorem
 is $1+1/d$. It turns out that,
actually,
$G_{\lambda}$ is summable to a higher power. The proof of this
is based on the parabolic version of
Gehring's lemma from \cite{GS_82}.

Introduce $\bQ$ as the set of cylinders $C_{R}(t,x)$, $R>0$,  
$t\geq0$, $x\in\bR^{d}$. For $Q=C_{R}(t,x)\in\bQ$ let
 $2Q=C_{2R}(t,x)$. If $Q\in\bQ$ and $Q=C_{R}(t,x)$,
we call $R$ the radius of $Q$.  

\begin{theorem}
                                           \label{theorem 9.3.2}
Let $ \lambda\in(0,\infty)$.
Then there exist  $d_{0}\in(1,d)$ and a constant    $N $, depending only
on $\delta,d,\uR,  p_{0}, \lambda$,  
such that for any     $Q\in\bQ$ of radius $R\leq\uR/2$ and $p\geq d_{0}+1 $, we have
\begin{equation}
                                          \label{10.14.1}
\| G_{\lambda}\|_{L_{p/(p-1)}(Q)}\leq N R^{-(d+2) /p }
\| G_{\lambda}\|_{L_{1}(2Q )} , 
\end{equation}
which is equivalently rewritten as
$$
\Big(\dashint_{Q}G^{p/(p-1)}_{\lambda}\,dxdt\Big)^{(p-1)/p}
\leq N\dashint_{2Q }G_{\lambda}\,dxdt.
$$
 
\end{theorem}

Proof. We basically follow the idea in \cite{FS_84}.
Take $Q\in\bQ$ of radius $R\leq\uR/2$ and
define 
recursively 
$$
\gamma^{1}=\inf\{t\geq0 :(t,x_{t})\in \bar Q\},\quad
\tau^{1}=\inf\{t\geq\gamma^{1} :(t,x_{t})\not\in 2Q\},
$$
$$
\gamma^{n+1}=\inf\{t\geq\tau^{n} :(t,x_{t})\in \bar Q\},\quad
\tau^{n+1}=\inf\{t\geq\gamma^{n+1} :(t,x_{t})\not\in 2Q\}.
$$
Then for any nonnegative Borel $f$ vanishing outside $Q$
with $\|f\|_{L_{d+1}(Q)}=1$
we have
$$ 
\int_{Q}f   G_{\lambda}(t,x)\,dxdt
=E\int_{0}^{\infty}e^{-\lambda t}f (t,x_{t})\,dt 
$$
$$
=\sum_{n=1}^{\infty}Ee^{- \lambda\gamma^{n}} 
E\Big(\int_{\gamma^{n}}^{\tau^{n}}e^{-\lambda(t-
 \gamma^{n} )}f(t,x_{t})\,dt
\mid \cF_{\gamma^{n}}\Big).
$$
Next we use the conditional  version of 
\eqref{9.29.2}   
 to see that the conditional expectations
above are less than $NR^{d/(d+1)}$. After that
we use the conditional  version of Corollary
\ref{corollary 7.29.1}  
to get that
$$
R^{2}\leq N E\Big(\int_{\gamma^{n}}^{\tau^{n}}
e^{-\lambda(t-\gamma^{n})} \,dt
\mid \cF_{\gamma^{n}}\Big).
$$
Then we obtain
$$
\int_{Q}f   G_{\lambda}(t,x)\,dxdt
\leq N  R ^{-(d+2)/(d+1)}
\sum_{n=1}^{\infty}E e^{-\lambda\gamma^{n}}
E\Big(\int_{\gamma^{n}}^{\tau^{n}}e^{-\lambda(t-\gamma^{n})} \,dt
\mid \cF_{\gamma^{n}}\Big)
$$
$$
=N  R ^{-(d+2)/(d+1)}
\sum_{n=1}^{\infty}E  
 \int_{\gamma^{n}}^{\tau^{n}}e^{-\lambda t } \,dt
$$
$$
\leq N  R ^{-(d+2)/(d+1)}
E \int_{0}^{\infty}e^{-\lambda t}I_{2Q}(t,x_{t}) \,dt
$$
$$
=N  R ^{-(d+2)/(d+1)}\int_{2Q}G_{\lambda}(t,x)\,dxdt.
$$

The arbitrariness of $f$ 
implies that
$$
\Big(\dashint_{Q}  G_{\lambda} ^{(d+1)/d}(t,x)\,dxdt
\Big)^{d/(d+1)}\leq N \dashint_{2Q}  G_{\lambda}(t,x)\,dxdt.
$$

Now the assertion of the theorem for $p=d_{0}$ follows directly from
the parabolic version of the famous Gehring's lemma
stated as Proposition 1.3  in \cite{GS_82}. 
For larger $p$ it suffices to use H\"older's inequality.
The theorem is proved.

The parameters $d_{0}$ and $N$ in Theorem
\ref{theorem 9.3.2} may depend
on $\delta$, $d$, $\uR$,  $p_{0}$,  
 $\lambda$.
What is important for the future is
that $d_{0}$ below is independent of $\lambda$.
\begin{theorem}
                                           \label{theorem 2.3.1}
There exists $d_{0}\in(1,d)$, depending only on
$\delta$, $d$, $\uR$,  $p_{0}$,  
 such that for any $p\geq d_{0}+1$ and
$\lambda>0$
$$
\int_{0}^{\infty}\int_{\bR^{d}}G_{\lambda}^{p/(p-1)}(t,x)
\,dxdt\leq N(\delta,d,\uR,p_{0},\lambda,p).
$$
Furthermore, the above constant
$N(\delta,d,\uR,p_{0},\lambda,p)$
can be taken in the form
$$
N(\delta,d,\uR,p_{0},p)\melambda_{p}^{(d+2)/(2p)-1},
$$
where
$$
\melambda_{p}=\lambda(1\wedge\lambda)^{d/(2p-d-2)}.
$$
\end{theorem}

Proof. Represent $\bR^{d+1}_{+}=[0,\infty)\times \bR$
as the union of countably many $Q_{1},Q_{2},...\subset \bQ$
of radius $\uR/2$
so that each point in $\bR^{d+1}_{+}$ belongs
to no more than $m(d)$  of the $2Q_{i}$'s and let $d_{0}$
be taken from Theorem \ref{theorem 9.3.2}
with   $\lambda=1$. Then for $\lambda\geq1$
$$
\|G_{\lambda}\|_{L_{p/(p-1)}(\bR^{d+1}_{+})}
\leq \|G_{1}\|_{L_{p/(p-1)}(\bR^{d+1}_{+})}
\leq \|\sum_{i}I_{Q_{i}}G_{ 1}\|_{L_{p/(p-1)}(\bR^{d+1}_{+})}
$$
$$
\leq\sum_{i} \| G_{1}\|_{L_{p/(p-1)}(Q_{i})}
\leq N\sum_{i} \| G_{1}\|_{L_{1}(2Q_{i})}\leq N_{1}
\| G_{1}\|_{L_{1}(\bR^{d+1}_{+})}=N_{1}.
$$

If $\lambda\in(0,1)$ we take nonnegative $f\in L_{p}(\bR^{d+1}_{+})$
and observe that
$$
J:=E\int_{0}^{\infty}e^{-\lambda t}f(t,x_{t})\,dt
=\sum_{n=0}^{\infty}e^{-\lambda n}
E\int_{n}^{n+1}e^{-\lambda (t-n)}f(t,x_{t})\,dt
$$
$$
\leq \sum_{n=0}^{\infty}e^{1-\lambda}e^{-\lambda n}
E\int_{n}^{n+1}e^{-  (t-n)}f(t,x_{t})\,dt.
$$
By the first case each expectation in the sum is dominated
by $N\|fI_{[n,n+1)}\|_{L_{p }}$. Therefore
$$
J\leq N\sum_{n=0}^{\infty} e^{-\lambda n}\|fI_{[n,n+1)}\|_{L_{p }}
\leq N(1-e^{-\lambda})^{-(p-1)/p}\|f\|_{L_{p}(\bR^{d}_{+})},
$$
where the second inequality follows
from H\"olders inequality.
This takes care of the case that $\lambda\in(0,1)$ 
in both statements of the theorem.

To prove the second statement in case $\lambda\geq1$
consider the process $(t,y_{t})$, where $y_{t}=\sqrt\lambda
x_{t/\lambda}$. We have
$$
y_{t}=\int_{0}^{t}\sigma_{s/\lambda}\,d(\sqrt\lambda w_{s/\lambda})
+\int_{0}^{t} \lambda^{-1/2} b_{s/\lambda}\,ds,
$$
where $\sqrt\lambda w_{s/\lambda}$ is a Wiener process and
$$
|\lambda^{-1/2} b_{s/\lambda}|\leq \lambda^{-1/2}h(s/
\lambda,x_{t/\lambda})=:\tilde h(s,y_{s}).
$$
Observe that
$$
\|\tilde h\|^{q_{0}}_{L_{p_{0},q_{0}}(C_{R}(t,x))}=
\sqrt\lambda\| h\|^{q_{0}}_{L_{p_{0},q_{0}}(C_{R/\sqrt\lambda}
(t\lambda,x/\sqrt\lambda))}\leq \bar b_{R/\sqrt\lambda}R
\leq \bar b_{R }R.
$$
It follows that the above theory is applicable to
$(t,y_{t})$ and provides estimates with the same constants
as for $(t,x_{t})$. In particular, for any $p\geq d_{0}+1$
(with $d_{0}$ found above) and Borel
nonnegative $f(t,x)$
$$
I:=E\int_{0}^{\infty}e^{-t}f(t,y_{t})\,dt
\leq NN(\delta,d,\uR,p_{0} ,p)\|f\|_{L_{p}
(\bR^{d+1}_{+})}.
$$
After that it only remains to note that
$$
I=\lambda E\int_{0}^{\infty}e^{-\lambda t}g(t,x_{t})\,dt,
$$
where $g(t,x)=f(\lambda t,\sqrt\lambda x)$ and
$$
\|f\|_{L_{p}
(\bR^{d+1}_{+})}=\lambda^{(d+2)/(2p)}\|g\|_{L_{p}
(\bR^{d+1}_{+})}.
$$
The theorem is proved.

Similar improvement of integrability occurs for
the Green's function of $x_{t}$ rather than $(t,x_{t})$.
Here is a straightforward consequence
of \eqref{10.4.2}.  

\begin{theorem}
                                       \label{theorem 10.4.1}
 Let $p\in[d,\infty)$. Then
for any $\lambda>0$ 
  there exists a nonnegative Borel function $g_{\lambda}( x)$
(Green's function of $ x_{\cdot} $)
on $\bR^{d }$
 such that   for any Borel nonnegative $f$ given on $\bR^{d }$
we have
$$
 E\int_{0}^{\infty}e^{-\lambda t}
f( x_{t}) \,d t =\int_{\bR^{d }}f( x)g_{\lambda}( x)\,dx ,
$$
\begin{equation}
                                                   \label{10.4.5}
 \|\Psi^{-d/p}_{\lambda} g_{\lambda}\|_{L_{p' }(\bR^{d}) } 
\le N\lambda^{-1+d/p}\nelambda^{-  d/(2p )},
\end{equation}
where   $\Psi _{\lambda}( x)=\exp(- \sqrt{\nelambda}
  |x| \theta/16)$,
  $p '=p /(p -1)$, and   $N$
depends only on $\delta$, $d$, $\uR,p, p_{0},\bar b_{\infty}$.
 
\end{theorem}

According to this theorem 
this Green's function is summable to the power $d/(d-1)$.
Again it turns out that this power can be increased.
If $B$ is an open ball in $\bR^{d}$ by $2B$ we denote
the concentric open ball of twice the radius of $B$.

\begin{theorem}
                                           \label{theorem 10.4.3}
Let $\lambda\in(0,\infty)$.
Then there exist  $d_{0}\in(1,d)$ and a constant    $N $, depending only
on $d,\delta$, $\uR$,    $\lambda$,
 such that for any   ball  $B$ of radius $R\leq \uR/2$ 
 and $p\geq d_{0}  $, we have
\begin{equation}
                                          \label{10.4.6}
\| g_{\lambda}\|_{L_{p/(p-1)}(B)}\leq N R^{-d /p }
\| g_{\lambda}\|_{L_{1}(2B )} , 
\end{equation}
which is equivalently rewritten as
$$
\Big(\dashint_{B}g^{p/(p-1)}_{\lambda}\,dx \Big)^{(p-1)/p}
\leq N\dashint_{2B }g_{\lambda}\,dx.
$$
 
\end{theorem}

Proof. We again follow the idea in \cite{FS_84}.
Take a ball  $B$ of radius $R\leq \uR/2$ and
define 
recursively 
$$
\gamma^{1}=\inf\{t\geq0 : x_{t} \in \bar B\},\quad
\tau^{1}=\inf\{t\geq\gamma^{1} : x_{t} \not\in 2B\},
$$
$$
\gamma^{n+1}=\inf\{t\geq\tau^{n} : x_{t} \in \bar B\},\quad
\tau^{n+1}=\inf\{t\geq\gamma^{n+1} : x_{t} \not\in 2B\}.
$$
Then for any nonnegative Borel $f$ vanishing outside $B$
with $\|f\|_{L_{d}(B)}=1$
we have
$$ 
\int_{B}f   g_{\lambda}(x)\,dx 
=E\int_{0}^{\infty}e^{-\lambda t}f ( x_{t})\,dt 
$$
$$
=\sum_{n=1}^{\infty}Ee^{-\lambda \gamma^{n}} 
E\Big(\int_{\gamma^{n}}^{\tau^{n}}e^{-\lambda(t-
 \gamma^{n} )}f( x_{t})\,dt
\mid \cF_{\gamma^{n}}\Big).
$$
Next we use the conditional  version of 
\eqref{9.29.5}
 to see that the conditional expectations
above are less than $NR $. After that
we use the conditional  version of Corollary
\ref{corollary 7.29.1} 
to get that   
$$
R^{2}\leq N E\Big(\int_{\gamma^{n}}^{\tau^{n}}
e^{-\lambda(t-\gamma^{n})} \,dt
\mid \cF_{\gamma^{n}}\Big).
$$
Then we obtain
$$
\int_{B}f   g_{\lambda}( x)\,dx 
\leq N  R ^{-1}
\sum_{n=1}^{\infty}E e^{-\lambda\gamma^{n}}
E\Big(\int_{\gamma^{n}}^{\tau^{n}}e^{-\lambda(t-\gamma^{n})} \,dt
\mid \cF_{\gamma^{n}}\Big)
$$
$$
=N  R ^{-1}
\sum_{n=1}^{\infty}E  
 \int_{\gamma^{n}}^{\tau^{n}}e^{-\lambda t } \,dt
$$
$$
\leq N  R ^{-1}
E \int_{0}^{\infty}e^{-\lambda t}I_{2B}( x_{t}) \,dt
=N  R ^{-1}\int_{2B}g_{\lambda}( x)\,dx .
$$

The arbitrariness of $f$ 
implies that
$$
\Big(\dashint_{B}  g_{\lambda} ^{d/(d-1)}(x)\,dx 
\Big)^{(d-1)/d}\leq N \dashint_{2B}  g_{\lambda}(x)\,dx,
$$
and again it only remains to use Gehring's lemma
in case $p=d$.
For larger $p$ it suffices to use H\"older's inequality.
The theorem is proved.

By mimicking the proof of Theorem \ref{theorem 2.3.1}
one gets its ``elliptic'' counterpart.

\begin{theorem}
                                      \label{theorem 2.3.2}
There exists $d_{0}\in(1,d)$, depending only on
$\delta$, $d$, $\uR$,  $p_{0}$,  
 such that for any $p\geq d_{0}$ and
$\lambda>0$
$$
 \int_{\bR^{d}}g_{\lambda}^{p/(p-1)}(x)
\,dx \leq N(\delta,d,\uR,p_{0} ,p)
\melambda_{e,p}^{d/(2p)-1},
$$
where
$$
\melambda_{e,p}=\lambda(1\wedge\lambda)^{d/(2p-d)}.
$$
\end{theorem}

\begin{remark}
                                           \label{remark 2.7.1}
Below by $d_{0}$ we denote the largest of the $d_{0}$'s
from Theorems \ref{theorem 2.3.1} and \ref{theorem 2.3.2}
and observe that, as the simple example of $a^{ij}=\delta^{ij}$
and $b\equiv0$
shows, $d_{0}>d/2$.

\end{remark}

Next, we present an improved mixed-norm parabolic Aleksandrov
estimates by following the
interpolation  arguments in Nazarov \cite{Na_15}.

\begin{lemma}
                                          \label{lemma 8.27.1}
Suppose that 
\begin{equation}
                                             \label{10.21.4}
p,q\in[1,\infty],\quad \frac{d_{0}}{p}+\frac{1}{q}=1.
\end{equation}
Then for any Borel $f(t,x)\geq0$  
\begin{equation}
                                             \label{8.27.1}
I:=E\int_{0}^{\infty}e^{-\lambda t}f(t,x_{t})\,dt
\leq N\melambda_{d_{0}+1}^{-(2d_{0}-d)/(2p)}\|f\|_{L_{p,q}},
\end{equation}
where $N=N(\delta,d,\uR,p ,p_{0} )$.
\end{lemma}

Proof. 
If $p=d_{0}+1$, then $q=d_{0}+1$ and 
\eqref{8.27.1} follows from  Theorem \ref{theorem 2.3.1}. 
In other terms, for any Borel $f(t,x)\geq0$  
$$
E\int_{0}^{\infty}e^{-\lambda t}f(t,x_{t})\,dt
=\int_{Q}G_{\lambda}(t,x)f(t,x)
\,dxdt
$$
$$
\leq N\melambda_{d_{0}+1}^{-(2d_{0}-d)/(2d_{0}+2)}
\|f\|_{L_{d_{0}+1} }.
$$

If $p=d_{0}$ and $q=\infty$ estimate \eqref{8.27.1}
follows from Theorem \ref{theorem 2.3.2} since
$$
I\leq E\int_{0}^{\infty}e^{-\lambda t}\sup_{s\geq0}f(s,x_{t})\,dt
=\int_{\bR^{d}}g_{\lambda}(x)\sup_{s\geq0}f(s,x)\,dx
$$
$$
\leq N\Big(\int_{\bR^{d}}\sup_{s\geq0}f^{d_{0}}(s,x)\,dx\Big)^{1/d_{0}}
=N\melambda_{e,d_{0}}^{-(2d_{0}-d)/(2d_{0})}\|f\|_{L_{d_{0},\infty}}
$$
and as is easy to check $\melambda_{e,d_{0}}=\melambda_{d_{0}+1}$.

If $p=\infty$ and $q=1$
 $$
I\leq \int_{0}^{\infty} \sup_{x}f(t,x)\,dt=\|f\|_{L_{\infty,1}}.
$$

We will use these facts in an interpolation argument.
In case $\infty>p>d_{0}+1$ we have $p>q$ and
set $\beta=p/(d_{0}+1)$ and $\alpha=\beta/(\beta-1)$.
Take a nonnegative $g(t)$ such that $ \big(f(t,x)g(t)
\big)/g(t)=f(t,x)$ ($0/0=0$) and use H\"older's inequality to conclude
that $I\leq I_{1}I_{2}$, where
$$
I_{1}=\Big(\int_{0}^{\infty} g^{-\alpha}(t)
\,dt\Big)^{1/\alpha},
$$
$$
I_{2}=\Big(E\int_{0}^{\infty}e^{-\lambda t}
g^{\beta}(t)f^{\beta}(t,x_{t})
\,dt\Big)^{1/\beta}
$$
$$
\leq N\melambda_{d_{0}+1}^{-(2d_{0}-d)/(2p)}\Big(\int_{0}^{\infty}g^{(d_{0}+1)\beta}(t)
\Big(\int_{\bR^d}f^{(d_{0}+1)\beta}(t,x)\,dx\Big)\,dt
\Big)^{1/(d_{0}\beta+\beta)}.
$$
For $g$ found from
$$
  g^{-\alpha}(t)=
g^{(d_{0}+1)\beta}(t)
\int_{\bR^d}f^{(d_{0}+1)\beta}(t,x)\,dx
$$
we get \eqref{8.27.1} and this takes care of the case
that $\infty>p>d_{0}+1$.

If $\infty>q>d_{0}+1$ we have $p<q$ and
set $\beta=q/(q-d_{0}-1)$ and $\alpha=\beta/(\beta-1)$.
Take a nonnegative $g(x)$ such that $ \big(f(t,x)g(x)
\big)/g(x)=f(t,x)$ ($0/0=0$) and use H\"older's 
inequality to conclude
that $I\leq I_{1}I_{2}$, where
$$
I_{1}=\Big(E\int_{0}^{\infty}e^{-\lambda
t}g^{-\beta}(x_{t})\,dt
\Big)^{1/\beta}
\leq N\melambda_{d_{0}+1}^{(d-2d_{0})/(2d_{0}\beta)}\Big(\int_{\bR^{d}}g^{-d_{0}\beta}(x)\,dx
\Big)^{1/(d_{0}\beta)},
$$
$$
I_{2}=\Big(E\int_{0}^{\infty}e^{-\lambda
t}g^{\alpha}(x_{t})
f^{\alpha}(t,x_{t})\,dt
\Big)^{1/\alpha}
$$
$$
\leq N\melambda_{d_{0}+1}^{(d-2d_{0})/(2d_{0}\alpha+2\alpha)}
\Big(\int_{\bR^{d}}g^{(d_{0}+1)\alpha}(x)
\Big(\int_{0}^{\infty}
f^{(d_{0}+1)\alpha}(t,x)\,dt\Big)dx
\Big)^{1/(\alpha d_{0}+\alpha )}.
$$
For $g$ found from
$$
g^{-d_{0}\beta}(x)=g^{(d_{0}+1)\alpha}(x)
 \int_{0}^{\infty}
f^{(d_{0}+1)\alpha}(t,x)\,dt 
$$
we get \eqref{8.27.1} after simple manipulations
and this proves the lemma.

Using this lemma instead of \eqref{9.29.2}
and just repeating the proof of Lemma \ref{lemma 8.22.1}
we come to a natural counterpart of the latter
and then by literally repeating the proof
of Theorem \ref{theorem 8.22.1} we come to the following   
result, that  is a version of Theorem 4.1 of Nazarov \cite{Na_15} in which
$d/p+1/q\leq 1$ that is stronger than ours,
but in which the assumption on $h$ is weaker.
A proper probabilistic 
version of Theorem 4.1 of Nazarov \cite{Na_15}
is found in \cite{Kr_20_2}. Recall that
$$
\bR^{d+1}_{+}=\{t\geq 0\}\cap\bR^{d+1}.
$$
\begin{theorem}
                                  \label{theorem 8.30.1}
Suppose   
\begin{equation}
                                                   \label{10.7.1}
p ,q \in [1,\infty],\quad
\nu:=1-\frac{d_{0}}{p }-\frac{1}{q }\geq 0.
\end{equation}
Then there is  $N=N(\delta,d,\uR,p,q,p_{0},\bar b_{\infty} )$ such that
for any $\lambda>0$ and   Borel nonnegative $f$  
we have 
\begin{equation}
                                                   \label{8.22.40}
 E\int_{0}^{\infty}e^{- \lambda t}
f(t,x_{t}) \,dt\leq 
N\melambda_{d_{0}+1}^{-\nu+(d-2d_{0})/(2p )}\|\Psi_{\lambda}^{1-\nu} f\|_
{L_{p ,q }(\bR^{d+1}_{+})},   
\end{equation}
where $\Psi _{\lambda}(t,x)=\exp(- 
\sqrt{\nelambda} (|x|+ \sqrt t)\bar\xi/16)$.
In particular, if $f$ is independent of $t$, $p\geq d_{0}$,
 and $q=\infty$
$$
E\int_{0}^{\infty}e^{- \lambda t}
f( x_{t}) \,dt\leq 
N\melambda_{d_{0}+1}^{-1+d/(2p)}\|\bar \Psi_{\lambda}^{d_{0}/p} f\|_
{L_{p   }(\bR^{d } )},
$$
where $\bar \Psi _{\lambda}( x)=\exp(- 
\sqrt{\nelambda}  |x| \bar\xi/16)$.  
\end{theorem}

\begin{theorem}
                                          \label{theorem 9.7.1}
Assume that \eqref{10.7.1} holds.
Then

(ii)
  for any
$n=1,2,...$, nonnegative Borel $f$ on $\bR^{d+1}_{+}$, and
 $T\leq 1$  we have
\begin{equation}
                                          \label{9.7.1}
E\Big[\int_{0}^{T}  
f(t,x_{t})\,dt\Big]^{n}\leq n!N^{n} 
T^{n\chi }\| \Psi^{(1-\nu)/n} _{1/T}
f\|^{n}_{L_{p,q}(\bR^{d+1}_{+}) },
\end{equation}
where $N=N(\delta,d,\uR,p,q,p_{0},\bar b_{\infty} )$ and 
$\chi=\nu+(2d_{0}-d)/(2p)$; 

(ii)   for any
  nonnegative Borel $f$ on $\bR^{d+1}_{+}$, and
 $T\geq 1$  we have
\begin{equation}
                                          \label{9.7.10}
I:=E \int_{0}^{T}  
f(t,x_{t})\,dt \leq N  T^{1-1/q}
 \| \Psi^{  1-\nu   } _{1}
f\| _{L_{p,q}(\bR^{d+1}_{+}) },
\end{equation}
where $N=N(\delta,d,\uR,p,q,p_{0},\bar b_{\infty} )$.
 
\end{theorem}

Proof. The proof of (i) proceeds by induction
on $n$ and is achieved by almost literally repeating
the proof of Theorem 2.7 of \cite{Kr_19_1}.
The induction hypothesis is that for all
$(t,x)\in \bR^{d+1}_{+}$ and $\kappa\in[0,1/n]$
\begin{equation}
                                                \label{10.28.1}
E\Big[\int_{0}^{T}  
f(t+s,x+x_{s})\,ds\Big]^{n}\leq n!N^{n} 
T^{n\chi } \Psi_{1/T}^{( \nu-1)\kappa n}
(t, x)\| \Psi^{( 1-\nu)\kappa} _{1/T}
f\|^{n}_{L_{p,q}(\bR^{d+1}_{+}) }.
\end{equation}

We will discuss in detail only the   case of $n=1$.
 In that case observe that
$\lambda:=1/T\geq 1 $
which allows us to use Theorem \ref{theorem 8.30.1}
in the following computations when $\kappa\in[0,1]$:
$$
 e^{-\lambda T}E \int_{0}^{T}  
f(t+s,x+x_{s})\,ds
\leq E \int_{0}^{\infty} e^{-\lambda s} 
f(t+s,x+x_{s})\,ds
$$
$$
\leq N \lambda^{-\chi}\|
f(t+\cdot,x+\cdot)
 \Psi^{(1-\nu)\kappa} _{\lambda}\|_{L_{p,q}(\bR^{d+1}_{+}) }
$$
$$
\leq N \lambda^{-\chi} \Psi^{(\nu-1)\kappa}_{\lambda}( x)\|
f 
 \Psi^{(1-\nu)\kappa} _{\lambda}\|_{L_{p,q}(\bR^{d+1}_{+}) },
$$
where the last inequality is due to the fact that,
for $x,y \in \bR^{d } $,
$\Psi  _{\lambda}(s, y)\leq \Psi  _{\lambda}(t+s, x+y)
\Psi ^{-1} _{\lambda}(t, x)$.
Since $\lambda T=1$, we get \eqref{10.28.1} with $n=1$. 

While proving \eqref{9.7.10} we may assume that  
  $T=k $, where $k\geq1$ is an integer. Then note that
owing to \eqref{9.7.1} for any integer $n\geq0$
$$
E \Big(\int_{n }^{ n+1  }  
f(t,x_{t})\,dt \mid\cF_{n }\Big)
\leq N\| \Psi_{1}^{1-\nu}f
I_{[n ,n+1 ]}
\|_{L_{p,q}}.
$$
Hence,
$$
I\leq N\sum_{n=0}^{k-1}
\| \Psi_{ 1}^{1-\nu}fI_{[n,n+1 ]}
\|_{L_{p,q}}=:NJ.
$$
If $p\geq q$, we have by H\"older's inequality
$$
J\leq k^{ (q-1)/q}\| \Psi_{1}^{1-\nu}
f \|_{L_{p,q}}.
$$
If $p<q$
$$
J\leq k^{ (p-1)/p}\Big(\int_{\bR^{d}}\sum_{n=0}^{k-1}
\Big(\int_{n }^{ n+1 }
 \Psi_{1}^{q(1-\nu)}
f^{q}(t,x)\,dt\Big)^{p/q}dx\Big)^{1/p}
 $$
$$
\leq k^{ (p-1)/p + (q-p)/(qp)}
\| \Psi_{1}^{1-\nu}
f \|_{L_{p,q}},
$$
which yields \eqref{9.7.10} since 
$(p-1)/p + (q-p)/(qp)=1-1/q$.
The theorem is 
proved.

Next theorem improves 
estimate \eqref{9.29.2} in what concerns
the restrictions on $p,q$.
\begin{theorem}
                                       \label{theorem 9.5.1}
Assume that \eqref{10.7.1} holds with $\nu=0$.
 Then 
for any $R\in(0,\bar R] $, $x$,  and Borel nonnegative $f$ 
given on $C_{R}$,
we have
\begin{equation}
                                                   \label{9.5.4}
E\int_{0}^{\tau_{R}(x)}f(t,x+x_{t})\,dt\leq
NR^{(2d_{0}-d)/p}\|f\|_{L_{p,q}(C_{R})},
\end{equation}
where $N=N(\delta,d,\uR,p, p_{0},\bar b_{\bar R},\bar R )$.

\end{theorem}

Proof.  Since $\tau_{R}(x)\leq R^{2}$, the left-hand side of
\eqref{9.5.4} is smaller than
$$
e^{\lambda R^{2}}E\int_{0}^{\infty}e^{-\lambda t}
I_{C_{R}}f(t,x+x_{t})\,dt
$$
for any $\lambda>0$. We estimate the last expectation
by using Theorem \ref{theorem 8.30.1}, observe
that, for $\lambda=R^{-2}$ we have $N\melambda_{d_{0}+1}
\geq \lambda$ owing to $R\leq \bar R$, and then
immediately come to \eqref{9.5.4}, however, with
$N$ depending on $\bar b_{\infty}$ in place of $\bar b_{\bar R}$.
 To see that this replacement is not needed, it suffices
to
note that the left-hand side of \eqref{9.5.4}
will not change if we replace $b_{t}$
with $b_{t}I_{t<\tau_{\bar R}(x)}$ which admits
the estimate
$$
|b_{t}I_{t<\tau_{\bar R}(x)}|\leq hI_{C_{\bar R}(0,-x)}(t,x_{t}).
$$
 The theorem
is proved.

Theorem \ref{theorem 9.7.1} allows us to prove It\^o's formula
for functions $u\in W^{1,2}_{p,q}(Q)$, where $Q$
is a domain in $\bR^{d+1}$ and 
$$
W^{1,2}_{p,q }(Q)=\{v: v, \partial_{t}v,
 Dv,  D^{2}v\in L_{p,q }(Q) \}
$$
with norm introduced in a natural way.
Before, the formula was known only for
(smooth, It\^o, and) $W^{1,2}_{d+1}$-functions
and processes with bounded drifts or
for $W^{2}_{d_{0}}$-functions in case the drift
of the process is dominated by $h(x_{t})$
with $h\in L_{d}$ (see \cite{Kr_19_1}).

The following extends Theorem 2.10.1  of \cite{Kr_77}.

\begin{theorem}
                                 \label{theorem 10.15.1}
Assume that \eqref{10.7.1} holds with $\nu=0$
and  $p<\infty$,
$q<\infty$. Let
$Q$ be a bounded domain in $\bR^{d+1}$, $0\in Q$,
$b$ be {\em bounded\/}, and  $u\in W^{1,2}_{p,q}(Q)\cap C(\bar Q)$. Then,
for $\tau$ defined as the first exit time of $(t,x_{t})$
from $Q$ and for all $t\geq0$,
$$
u(t\wedge\tau,x_{t\wedge\tau})
=u(0,0)+\int_{0}^{t\wedge\tau}D_{i}u(s,x_{s})\,dm^{i}_{s}
$$
\begin{equation}
                                      \label{10.15.1}
+\int_{0}^{t\wedge\tau}[
\partial_{t}u(s,x_{s})+ a^{ij}_{s}D_{ij}u(s,x_{s})
+b^{i}_{s}D_{i}u(s,x_{s})]\,ds
\end{equation}
and the stochastic integral above is a square-integrable
martingale.
\end{theorem}

Proof. First assume that $u$ is smooth and its derivatives
 are bounded. Then  
\eqref{10.15.1} holds by It\^o's formula and, moreover,
by denoting $\tau^{n}= n\wedge\tau$
for any $n\geq0$ we have
$$
E\int_{ \tau^{n}}^{\tau^{n+1}}|Du(s,x_{s})|^{2}\,ds
\leq NE\Big(\int_{\tau^{n}}^{ \tau^{n+1}
}D_{i}u(s,x_{s})\,dm^{i}_{s}\Big)^{2}
$$
$$
=NE\Big(u( \tau^{n+1},x_{ \tau^{n+1}})-
u( \tau^{n },x_{ \tau^{n }})
$$
$$-
\int_{ \tau^{n}}^{\tau^{n+1}}[
\partial_{t}u(s,x_{s})+ a^{ij}_{s}D_{ij}u(s,x_{s})
+b^{i}_{s}D_{i}u(s,x_{s})]\,ds\Big)^{2}
$$
$$
\leq N\sup_{\bar Q}|u|
+NE\Big(\int_{ \tau^{n}}^{\tau^{n+1}}I_{Q}
\big(|\partial_{t}u|+|Du|+|D^{2}u|\big)(s,x_{s})\,ds\Big)^{2}.
$$
Since $Q$ is bounded, $\tau$ is bounded as well
and
in light of Theorem \ref{theorem 9.7.1} we conclude that
\begin{equation}
                                      \label{10.15.2}
E\int_{0}^{\tau}|Du(s,x_{s})|^{2}\,ds\leq N
\sup_{\bar Q}|u|+N \| \partial_{t}u,
 Du, D^{2}u\|_{L_{p,q}(Q)},
\end{equation}
where $N$ are independent of $u$ and $Q$
as long as the size of $Q$ in the $t$-direction
is under control.
Owing to Fatou's theorem,
this estimate is also true for those $u\in W^{1,2}_{p,q}(Q)
\cap C(\bar Q)$
that can be approximated uniformly and in the 
$W^{1,2}_{p,q}(Q)$-norm by smooth functions with bounded
derivatives (recall that $p<\infty$,
$q<\infty$). For our $u$
there is no guarantee that such approximation is possible.
However, mollifiers do such approximations
 in any subdomain $Q'\subset \bar Q'\subset
Q$. Hence, \eqref{10.15.2} holds for our $u$ if we replace $Q$
by $Q'$
(containing $(0,0)$). Setting $Q' \uparrow Q$  proves \eqref{10.15.2}
in the generals case and proves the last assertion
of the theorem.

After that \eqref{10.15.1} with $Q'$ in place of $Q$
is proved by routine approximation of $u$ by smooth 
functions. Setting $Q' \uparrow Q$ finally proves \eqref{10.15.1}.
The theorem is proved.

\mysection{Application to parabolic equations}
                                             \label{section 10.13.1}

Fix a constant $\delta\in(0,1)$ and recall that by
$\bS_{\delta}$ we denote the set of $d\times d$-symmetric
matrices whose eigenvalues are between
$\delta$ and $\delta^{-1}$. In this section
we impose the following.
\begin{assumption}
                                         \label{assumption 10.12.1}

(i) On $\bR^{d+1}$
 we are given Borel measurable $\sigma(t,x)$
and $b(t,x)$ with values  in $\bS_{\delta}$
and in $\bR^{d}$ respectively.

(ii) We are given $p_{0},q_{0} \in[1,\infty)$  
satisfying \eqref{5.10.1} and a function $h(t,x)$
satisfying \eqref{8.19.1} 
and such that $|b|\leq h$.

(iii) Assumption \ref{assumption 12.18.2}
is satisfied.
\end{assumption}

Introduce $a=\sigma^{2}$ and set
$$
Lu(t,x)=(1/2)a^{ij}(t,x)D_{ij}u(t,x)+b^{i}(t,x)D_{i}u(t,x).
$$

For a domain $Q\subset\bR^{d+1}$
one denotes by $\partial'Q$ its parabolic boundary
defined as the set of all points on $\partial Q$
which are endpoints of continuous curves of type $(t,x_{t})$,
$t\in[a,b]$, which
start in $Q$ and belong to $Q$ for all $t<b$.

The following has the same flavor as Nazarov's
Theorem 4.1 of  \cite{Na_15}
or Theorem 4.3 of  \cite{Kr_20_2}. We get wider
range of $p,q$ on the account
of restricting $b$.
Here is a qualitative form
of the maximum principle.

\begin{theorem}
                                        \label{theorem 10.14.1}
Let $0<R\leq\bar R$, domain $Q\subset C_{R}$,  and
assume that \eqref{10.7.1} holds with $\nu=0$, $p<\infty$,
$q<\infty$, and that
we are given a function 
$u\in W^{1,2}_{p,q,\loc}(Q)\cap C(\bar Q)$.
  Take a function
$c\geq 0$ on $Q$. Then on $ Q$
\begin{equation}
                                               \label{10.14.10}
u \leq NR^{(2d_{0}-d)/p}
\|I_{Q,u>0}(\partial_{t}u+Lu-cu)_{-}\|_{L_{p,q} }
+\sup_{\partial'Q}u_{+},
\end{equation}
where $N=N(\delta,d,\uR,\bar R,p, p_{0},\bar b_{\bar R})$.
In particular (the maximum principle),
if $\partial_{t}u+Lu-cu \geq0$ in $Q$ and $u\leq0$
on $\partial'Q$, then $u\leq 0$ in $Q$.
\end{theorem}

Proof. Obviously the right-hand side of
\eqref{10.14.10} decreases if we replace $c$ with zero.
Hence we may assume that $c=0$. 
Also, we need to prove \eqref{10.14.10}
only in $Q\cap\{u>0\}$ on the parabolic boundary of
which either $u=0$ or $u\leq \sup_{\partial'Q}u_{+}$.
Therefore, we may assume that $u>0$ in $Q$.

Then for $\varepsilon>0$ define $Q^{\varepsilon}$
as the collection of $(t,x)\in Q$ such that
the closed ball in $\bR^{d+1}$ centered at $(t,x)$
with radius $\varepsilon$ lies in $Q$.
Obviously $Q^{\varepsilon}$ is   open.
It is not hard to prove (see, for instance,
the proof of Lemma 3.1.13 in \cite{Kr_18})
that $\dist(\partial' Q,\partial'Q^{\varepsilon})=\varepsilon$.
It follows, owing to the continuity of $u$
and the monotone convergence theorem, that
it suffices to prove \eqref{10.14.10}
with $Q^{\varepsilon}$ in place of $Q$.
As a consequence of that we may assume that
$u\in W^{1,2}_{p,q}(Q)$.

This  gives us the opportunity
to replace $L$ in  \eqref{10.14.10}
with $L_{n}:=I_{|b|\geq n}\Delta+I_{|b|< n}L$
and then pass to the limit
by the dominated convergence and monotone convergence theorems.
Hence, we may assume that $b$ is bounded.
In this situation for fixed $(t,x)\in Q$ by It\^o's formula
we have
\begin{equation}
                                                       \label{2.8.1}
u(t,x)=Eu(t+\tau,x_{\tau})+E\int_{0}^{\tau}f(t+s,x_{s})\,ds,
\end{equation}
where $f=-(\partial_{t}u+Lu)$, $x_{s}$ is a 
solution of \eqref{11.29.2}
on a probability space, and $\tau$ is the first exit time
of $(t+s,x_{s})$ from $Q$.
  After that, to prove \eqref{10.14.10},
 it only remains
to use Theorem \ref{theorem 9.5.1} and the fact that %%%
$(t+\tau
,x_{\tau})\in 
\partial'Q$.
The theorem is proved.

\begin{remark}
The statement of Theorem \ref{theorem 10.14.1}
may look futile because there is no guarantee
that even for an $u\in W^{1,2}_{p,q}(Q)$ the norm
in \eqref{10.14.10} is finite. However, there is an important case
(see \cite{Kr_21}) when $a=(\delta^{ij})$, $b\in L_{d+1}$,
and the equation $\partial_{t}u+\Delta u+b^{i}D_{i}u=f$
has solutions in $W^{1,2}_{p}$ with $Du\in L_{r}$,
$p\in (1,d+1)$, $r=p(d+1)/(d+1-p)$,
on the account of $f\in L_{r}$.
In that case $b^{i}D_{i}u\in L_{p}$.
\end{remark}

The last result  we present is needed
for constructing the theory
of fully nonlinear parabolic equations
with drift terms in $L_{p,q}$ following the path
in \cite{Kr_18}.    

\begin{theorem}
                                   \label{theorem 10.14.2}  
Assume that \eqref{10.7.1} holds with $p<\infty$,
$q<\infty$
and take $R\in(0,\infty]$. Then there exists
  constants $N,\kappa>0$, depending only on
$ d,\delta,\uR, p,q, p_{0},\bar b_{\infty} $,
such that for any $\lambda\geq 1 $
and $u\in W^{1,2}_{p,q,\loc}(C_{R})\cap C(\bar C_{R})$ 
($C_{\infty}=\bR^{d+1}_{+}$, $C(\bR^{d+1}_{+})$
is the set of bounded continuous functions on $\bR^{d+1}_{+}$)
 we have
$$
\lambda\|u_{+}\|_{L_{p,q}(C_{R/2})}
\leq 
 N\| (\lambda u-Lu-\partial_{t}u)_{+}\|_{L_{p,q}(C_{R}) } 
$$
\begin{equation}
                                           \label{10.14.40}  
+ N\lambda R^{d/p+2/q}e^{-\kappa R\sqrt{\lambda} }
\sup_{\partial' C_{R}}u_{+} ,
\end{equation}
 where the last term should be dropped if $R=\infty$.
\end{theorem}

Proof.   By having in mind 
the possibility to approximate
$C_{R}$ from inside by similar domains,
we see that we may assume that 
$u\in W^{1, 2}_{p,q }(C_{R})$ and $R<\infty$.
Then as in the proof of Theorem \ref{theorem 10.14.1}
we reduce the general case to the one in which
$b$ is bounded and after that we may assume that $u$ is smooth.
In that case, for $(t,x)\in C_{R/2}$, 
similarly to  \eqref{2.8.1}
$$
u(t,x)=E e^{-\lambda\tau_{R}}u(t+\tau_{R},x_{\tau_{R}})
-E \int_{0}^{\tau_{R}}e^{-\lambda t}f(t+s,x_{s})\,ds
$$
$$
=:I(t,x)+J(t,x),
$$
where $f=\lambda u-Lu-\partial_{t}u$,
$\tau_{R}$ is the first exit time of $(t+s,x_{s})$
from $C_{R}$ and $x_{s}$ is a solution of \eqref{11.29.2}. 

Here, thanks to \eqref{8.20.1}  
\begin{equation}
                                               \label{10.4.30}
I(t,x)\leq Ne^{-\kappa  R\sqrt{\lambda} }\sup_{\partial'C_{R}}u_{+},
\quad\|I_{C_{R/2}}I_{+}\|_{L_{p,q}  }\leq
NR^{d/p+2/q}e^{-\kappa  R
\sqrt{\lambda} }\sup_{\partial'C_{R}}u_{+},
\end{equation}
where $N,\kappa>0$ depend only on    
$ d,\delta,\uR, p,q, p_{0}  $.

To estimate $J$ we  define $f$ as zero outside $C_{R}$
and observe that
$$
J(t,x)\leq E \int_{0}^{\infty}e^{-\lambda t}f_{+}(t+s,x_{s})\,ds=:
\bar J(t,x)
$$
By Theorem \ref{theorem 8.30.1} 
we have
\begin{equation}
                                          \label{2.8.3}
\bar J(t,x)\leq 
N \lambda^{-\eta}
\|\Psi_{\lambda}^{1-\nu} f_{+}(t+\cdot,x+\cdot)\|_{L_{p ,q }
(\bR^{d+1}_{+})},
\end{equation}
where $\eta= \nu+( 2d_{0}-d)/(2p )$.

If $p\geq q$,
\eqref{2.8.3} implies that
\begin{equation}
                                         \label{11.3.10}
\int_{\bR^{d}}|\lambda^{\eta}\bar J(t,x)|^{p}\,dx
\leq N 
\int_{\bR^{d}}\Big(\int_{0}^{\infty}
F^{q/p}(t,s,x)\,ds\Big)^{p/q}dx ,
\end{equation}
where
$$
F(t,s,x)=\int_{\bR^{d}}\Psi^{(1-\nu)p}_{\lambda}(s,y)f^{p} 
(t+s,x+y)\,dy.
$$
By Minkowski's inequality the   integral
on the right in \eqref{11.3.10} is dominated by
$$
\Big(\int_{0}^{\infty}\Big(\int_{\bR^{d}}F(t,s,x)\,dx\Big)^{q/p}
\,ds\Big)^{p/q},
$$
where
$$
\int_{\bR^{d}}F(t,s,x)\,dx=\int_{\bR^{d}}f^{p} (t+s,y)\,dy
\int_{\bR^{d}}\Psi^{(1-\nu)p}_{\lambda}(s,y)\,dy
$$
$$
\leq N \lambda^{-d/2}e^{-\mu \sqrt{\lambda s}}
\int_{\bR^{d}}f^{p} (t+s,y)\,dy,
$$
with $\mu=\mu(\delta,p,q,\uR) >0$. Below by $\mu$
we denote all such constants.
It follows that
$$
\int_{\bR^{d}}|\lambda_{d_{0}+1}^{\eta}\bar J(t,x)|^{p}\,dx
\leq N \lambda^{-d/2}\Big(\int_{0}^{\infty}
e^{-\mu \sqrt{\lambda s}}\Big(
\int_{\bR^{d}}f^{p}(t+s,y)\,dy\Big)^{q/p}ds\Big)^{p/q},
$$
$$
\|\lambda^{\eta}\bar J\|^{q}_{L_{p,q}(\bR^{d+1}_{+})}
\leq N \lambda^{-1-qd/(2p)}\|f \|^{q}_{L_{p,q}(\bR^{d+1}_{+})},
$$
which along with \eqref{10.4.30} yield \eqref{10.14.40}.

If $q\geq p$,
$$
\int_{0}^{\infty}|\lambda^{\eta}\bar J(t,x)|^{q}\,dt
\leq N\int_{0}^{\infty}\Big(\int_{\bR^{d}}
 F^{p/q}\,(t,x,y) dy\Big)^{q/p}dt
$$
where
$$
F(t,x,y)=\int_{0}^{\infty}
\Psi^{(1-\nu)q}_{\lambda}(s,y)f^{q} 
(t+s,x+y)\,ds.
$$
By Minkowski's inequality
$$
\Big(\int_{0}^{\infty}|\lambda^{\eta}
\bar J(t,x)|^{q}\,dt\Big)^{p/q}
\leq N\int_{\bR^{d}}\Big(\int_{0}^{\infty}
F(t,x,y)\,dt\Big)^{p/q}dy,
$$
where
$$
\int_{0}^{\infty}
F(t,x,y)\,dt\leq
\int_{0}^{\infty}
 f^{q} 
( s,x+y)\,ds\int_{0}^{\infty}
\Psi^{(1-\nu)q}_{\lambda}(s,y) \,ds
$$
$$
\leq N\lambda^{-1}e^{-\mu\sqrt\lambda|y|}\int_{0}^{\infty}
 f^{q} 
( s,x+y)\,ds.
$$
Hence,
$$
\Big(\int_{0}^{\infty}|\lambda^{\eta}\bar J(t,x)|^{q}\,dt\Big)^{p/q}
\leq N\lambda^{-p/q}\int_{\bR^{d}}e^{-\mu\sqrt\lambda |y|}\Big(
\int_{0}^{\infty}
 f^{q} 
( s,x+y)\,ds\Big)^{p/q}\,dy,
$$
$$
\|\lambda^{\eta}\bar J\|_{L_{p,q}(\bR^{d+1}_{+})}^{p}
\leq N\lambda^{-p/q-d/2}\|f \|_{L_{p,q}(\bR^{d+1}_{+})}^{p}
$$
and we again come to \eqref{10.14.40}.
The theorem is proved.

The full strength of Theorem \ref{theorem 10.14.2}
is seen in the theory of fully nonlinear equations.
But even for linear ones one gets a nontrivial information
as, for instance, in the following theorem which, in particular,
implies that
the operator $L+\partial_{t}$ with the domain
$$
\{u\in W^{1,2}_{p,q,\loc}(C_{R})
\cap C(\bar C_{R}):Lu+\partial_{t}u\in L_{p,q}(C_{R}),
u_{\big|\partial'C_{R} }=0\}
$$
is a closed operator in $L_{p,q}(C_{R})$.

\begin{theorem}
                                               \label{theorem 2.8.1}
Assume that \eqref{10.7.1} holds with $p<\infty$,
$q<\infty$
and take $R\in(0,\infty)$. Suppose we are given
 $u_{0},u_{1},...\in W^{1,2}_{p,q,\loc}(C_{R})
\cap C(\bar C_{R})$ and $f\in L_{p,q}(C_{R})$ such that
$f_{n}:=Lu_{n}+\partial_{t}u_{n}\in L_{p,q}(C_{R})$ for $n\geq0$,
$$
\sup_{n\geq1}\sup_{\partial'C_{R}}|u_{n}|<\infty,
\quad \|f_{n}-f\|_{L_{p,q}(C_{R})}+\|u_{n}-u_{0}\|_{L_{p,q}(C_{R})}
\to 0
$$
as $n\to\infty$. Then $Lu_{0}+\partial u_{0}=f$ in $C_{R}$.
\end{theorem}

Proof. Take a smooth $\psi$ on $C_{R}$ and apply
\eqref{10.14.40} to $u_{n}-u+\psi/\lambda$
in place of $u$. Then pass to the limit as $n\to \infty$
to find
$$
\|\psi_{+}\|_{L_{p,q}(C_{R/2})}
\leq N_{1}\lambda R^{d/p+2/q}
e^{-\kappa R\sqrt{\lambda} }
$$
$$
+N_{2}\|\psi-f+Lu_{0}+\partial_{t}u_{0}
-(L\psi+\partial_{t}\psi)/\lambda
\|_{L_{p,q}(C_{R})},
$$
where $N_{2}$ is independent of $\lambda$ and $\psi$
and $N_{1}$ is independent of $\lambda$. By setting
$\lambda\to\infty$ we get
$$
\|\psi_{+}\|_{L_{p,q}(C_{R/2})}
\leq N_{2}\|\psi-f+Lu_{0}+\partial_{t}u_{0}\|_{L_{p,q}(C_{R})}.
$$
This is true if $\psi$ is smooth enough and by approximation
is true for any $\psi\in L_{p,q}(C_{R})$. For
$\psi=f-Lu_{0}-\partial_{t}u_{0}$ we get that
 $f-Lu_{0}-\partial_{t}u_{0}\leq 0$
in $C_{R/2}$. The reader understands that here as well as in  
\eqref{10.14.40} one can take any number $<R$ in place of $R/2$.
Hence, $f-Lu_{0}-\partial_{t}u_{0}\leq 0$ in $C_{R}$. Passing to
$-u_{n}$, $-f$ yields $f-Lu_{0}-\partial_{t}u_{0}\geq 0$
and proves the theorem.

\end{document}